\documentclass[12pt,twoside]{article}

\usepackage{a4,graphicx,amsmath,amsfonts,amssymb,mathrsfs,bbm}

\newcommand{\real}{\mathbbm{R}}
\newcommand{\complex}{\mathbbm{C}}

\newcommand{\ltwo}{\mathscr{L}_2[0,\infty)}
\newcommand{\htwo}{\mathscr{H}_2}
\newcommand{\hinf}{\mathscr{H}_{\infty}}
\newcommand{\htop}{\rm H}  
\renewcommand{\top}{{\rm T}}

\newtheorem{definition}{Definition}
\newtheorem{theorem}{Theorem}

\newtheorem{lemma}{Lemma}

\begin{document}

\thispagestyle{plain}

\begin{center}
  {\bf \Large Frequency domain integrals for \\[0.5ex]
    stability preservation in Galerkin-type \\[1.5ex]
    projection-based model order reduction}

\vspace{10mm}

{\large Roland~Pulch} \\[1ex]
{\small Institute of Mathematics and Computer Science, 
University of Greifswald, \\
Walther-Rathenau-Str.~47, 17489 Greifswald, Germany. \\
Email: {\tt roland.pulch@uni-greifswald.de}}

\end{center}

\bigskip\bigskip

%%%%%%%%%%%%%%%%%%%%%%%%%%%%%%%%%%%%%%%%%%%%%%%%%%%%%%%%%%%%%%%%%%%%%%%%%%%%%
%%%                         Abstract                                      %%%
%%%%%%%%%%%%%%%%%%%%%%%%%%%%%%%%%%%%%%%%%%%%%%%%%%%%%%%%%%%%%%%%%%%%%%%%%%%%%

\begin{center}
{Abstract}

\begin{tabular}{p{13cm}}
We investigate linear dynamical systems consisting of
ordinary differential equations with high dimensionality.
Model order reduction yields alternative systems of much lower dimensions.
However, a reduced system may be unstable, although
the original system is asymptotically stable.
We consider projection-based model order reduction of Galerkin-type.
A transformation of the original system ensures that any reduced system
is asymptotically stable.
This transformation requires the solution of a high-dimensional
Lyapunov inequality.
We solve this problem using a specific Lyapunov equation.
Its solution can be represented as a matrix-valued integral
in the frequency domain.
Consequently, quadrature rules yield numerical approximations,
where large sparse linear systems of algebraic equations
have to be solved.
We analyse this approach and show a sufficient condition on the error 
to meet the Lyapunov inequality.
Furthermore, this technique is extended to systems
of differential-algebraic equations with strictly proper transfer
functions by a regularisation.
Finally, we present results of numerical computations for
high-dimensional examples, which indicate the efficiency of
this stability-preserving method.

\bigskip


\bigskip

% 65L05 NA, ODEs, initial value problems
% 65L20 NA, ODEs, stability and convergence of numerical methods
% 65F10 NA, Numerical linear algebra, iterative methods for linear systems
% 34C20 ODEs, transformation and reduction of equations and systems
% 34D20 ODEs, stability theory, stability
% 93B40 Systems theory, control, computational methods
% 93D20 Systems theory, control, asymptotic stability
\end{tabular}
\end{center}

\clearpage

\markboth{\em R.~Pulch}{\em Frequency Domain Integrals for Stability Preservation in MOR}

%%%%%%%%%%%%%%%%%%%%%%%%%%%%%%%%%%%%%%%%%%%%%%%%%%%%%%%%%%%%%%%%%%%%%%%%%%%%%
%%%                        Introduction                                   %%%
%%%%%%%%%%%%%%%%%%%%%%%%%%%%%%%%%%%%%%%%%%%%%%%%%%%%%%%%%%%%%%%%%%%%%%%%%%%%%

\section{Introduction}
Numerical simulation represents a main tool for the investigation
of dynamical systems in science, engineering and other fields
of application.
High-fidelity modelling is required to obtain detailed information
on complex problems.
However, the high-fidelity systems may have a huge number of
state variables, which makes a numerical simulation expensive or
even infeasible. 
Hence methods of model order reduction (MOR) are applied to decrease
the dimensionality of the dynamical systems,
see~\cite{antoulas,benner-mehrmann,schilders}.
Yet the reduced system has to reproduce the quantities of interest
sufficiently accurate.

We consider linear systems of ordinary differential equations (ODEs),
which are asymptotically stable. 
Projection-based MOR determines linear ODEs with a lower dimensionality.
However, the reduced system may be unstable and thus useless.
Some solutions become unbounded for unstable systems in the time domain. 
Furthermore, error bounds, which follow from the transfer functions
in the frequency domain, are not available any more.
Hence stability-preserving MOR methods are essential to generate
appropriate reduced systems.

The balanced truncation technique, see~\cite{gugercin-antoulas},
always produces stable reduced systems, while the computational effort
is often relatively large.
Krylov subspace techniques, see~\cite{freund}, are less expensive,
whereas stability can easily be lost.
A stability preservation of a Krylov subspace approach is achieved by
special assumptions and methods in~\cite{ionescu}.
A post-processing, which works on the poles of the transfer function,
can recondition the stability, see~\cite{bai-freund}.

We employ projection-based MOR of Galerkin-type, where each scheme
is defined by a single orthogonal projection matrix.
Important methods are the one-sided Arnoldi algorithm and
the proper orthogonal decomposition (POD), for example.
% Comment: Freund and Antoulas are cited later in the paper !
A transformation of the system of ODEs guarantees the stability
of any reduced system, see~\cite{castane-selga,prajna,pulch-arxiv}.
This technique was also applied to a stochastic Galerkin projection
in~\cite{pulch-augustin}, where the stability of larger systems
than the original ODEs is ensured. 
In our MOR methods,
the main effort consists in solving a single high-dimensional
Lyapunov inequality,
where the efficient numerical solution is critical.

The Lyapunov inequality can be satisfied by the approximate solution of a
high-dimensional Lyapunov equation.
Therein, we perform a simple but effective choice of an input matrix.
We prove an error bound on the approximation, which is sufficient
for achieving the Lyapunov inequality.
The solution of the Lyapunov equation also represents a
matrix-valued integral in the frequency domain.
Phillips and Silveira~\cite{phillips-silveira} computed integrals of
this type approximately by a quadrature rule.
Quadrature methods were also considered for such integrals
in~\cite{benner-schneider,breiten}.
We use this approach to construct a stability-preserving MOR technique.
Therein, not the solution of the Lyapunov equation itself is required
but an associated matrix-matrix product with a small number of columns.
Now large sparse linear systems of algebraic equations have to be solved,
where the linear dynamical system yields the coefficient matrices.

Furthermore, we extend this stability-preserving MOR method to
systems of differential-algebraic equations (DAEs).
M\"uller~\cite{mueller} investigated a regularisation technique,
which changes an asymptotically stable DAE system into an
asymptotically stable ODE system.
Now the stability-preserving MOR applies to this system of ODEs.
In the case of DAEs with a strictly proper transfer function, 
we show that the additional regularisation error converges to zero
in dependence on a regularisation parameter.

The paper is organised as follows.
We introduce the considered MOR methods in Section~\ref{sec:problem-def}. 
The stability-preserving transformation and the Lyapunov equations
are discussed in Section~\ref{sec:preservation}.
We arrange the frequency domain integrals and analyse the usage
of quadrature rules.
In Section~\ref{sec:daes}, the stability-preserving approach is
transferred to systems of DAEs.
Finally, we present results of numerical experiments in
Section~\ref{sec:examples}, where an ODE system and a DAE system
are examined. 

%%%%%%%%%%%%%%%%%%%%%%%%%%%%%%%%%%%%%%%%%%%%%%%%%%%%%%%%%%%%%%%%%%%%%%%%%%%%%
%%%                     Problem Definition                                %%%
%%%%%%%%%%%%%%%%%%%%%%%%%%%%%%%%%%%%%%%%%%%%%%%%%%%%%%%%%%%%%%%%%%%%%%%%%%%%%

\section{Model order reduction and stability}
\label{sec:problem-def}
Projection-based MOR of linear dynamical systems is closely related to
stability properties, which are reviewed in this section.

%%%%%%%%%%%%%%%%%%%%%%%%%%%%%%%%%%%%%%%%%%%%%%%%%%%%%%%%%%%%%%%%%%%%%%%%%%%%%
\subsection{Linear dynamical systems and stability}
We consider linear time-invariant systems in the form 
\begin{equation} \label{linear-system}
  \begin{array}{rcl}
    E \dot{x}(t) & = & A x(t) + B u(t) \\
    y(t) & = & C x(t) \\
  \end{array}
\end{equation}
with state/inner variables $x : [0,\infty) \rightarrow \real^n$,
inputs $u : [0,\infty) \rightarrow \real^{n_{\rm in}}$ and
outputs $y : [0,\infty) \rightarrow \real^{n_{\rm out}}$.
The system includes constant matrices $A,E \in \real^{n \times n}$,
$B \in \real^{n \times n_{\rm in}}$ and $C \in \real^{n_{\rm out} \times n}$.
If the mass matrix~$E$ is non-singular, then the system~(\ref{linear-system})
consists of ordinary differential equations (ODEs).
If the mass matrix~$E$ is singular, then differential-algebraic
equations (DAEs) are given.
The pair $(E,A)$ is called a matrix pencil. 
We assume that the matrix pencil is regular, i.e.,
$\det ( \lambda E - A ) \neq 0$ for some $\lambda \in \complex$.
ODEs always yield a regular matrix pencil.
We add predetermined initial values $x(0)=x_0$,
which are assumed to be consistent in the case of DAEs. 

In the frequency domain, a transfer function describes
the input-output behaviour of the system~(\ref{linear-system}) completely,
see~\cite{antoulas}.
This transfer function
$H : \complex \backslash \Sigma \rightarrow \complex^{n_{\rm out} \times n_{\rm in}}$
reads as
\begin{equation} \label{transfer}
  H(s) = C ( s E - A )^{-1} B
  \qquad \mbox{for} \;\; s \in \complex \backslash \Sigma .
\end{equation}
The mapping~(\ref{transfer}) is a rational function with a
finite set of poles~$\Sigma \subset \complex$.
The magnitude of a transfer function can be characterised by norms
in Hardy spaces.
The $\htwo$-norm is defined by, see~\cite[p.~92]{shmaliy},
\begin{equation} \label{h2-norm}
  \left\| H \right\|_{\htwo} =
  \sqrt{ \frac{1}{2\pi} \int_{-\infty}^{+\infty}
    \left\| H({\rm i}\omega) \right\|_{\rm F}^2 \; {\rm d}\omega } 
  %\sqrt{ \int_{-\infty}^{+\infty} {\rm trace}
  %  ( H({\rm i}\omega)^* H({\rm i}\omega) ) \; {\rm d}\omega } \; .
\end{equation}
with ${\rm i} = \sqrt{-1}$ and the Frobenius (matrix) norm
$\| \cdot \|_{\rm F}$ provided that the integral exists.

The stability issues of the system~(\ref{linear-system}) are independent of
the definition of inputs or outputs. 
To discuss the stability, we recall some general properties of matrices
in the following definitions.

%%%%%%%%%%%%%%%%%%%%%%%%%%%%%%%%%%%%%%%%%%%%%%%%%%%%%%%%%%%%%%%%%%%%%%%%%%%%%
\begin{definition} \label{def:spectral}
  Let $A \in \real^{n \times n}$ and $\lambda_1,\ldots,\lambda_n \in \complex$
  be its eigenvalues.
  The {\em spectral abscissa} of the matrix~$A$ is the real number
  $$ \alpha (A) = \max
  \left\{ {\rm Re}(\lambda_1) , \ldots , {\rm Re} (\lambda_n) \right\} . $$
\end{definition}
%%%%%%%%%%%%%%%%%%%%%%%%%%%%%%%%%%%%%%%%%%%%%%%%%%%%%%%%%%%%%%%%%%%%%%%%%%%%%

%%%%%%%%%%%%%%%%%%%%%%%%%%%%%%%%%%%%%%%%%%%%%%%%%%%%%%%%%%%%%%%%%%%%%%%%%%%%%
%\begin{definition} \label{def:stable-matrix}
%  A matrix $A \in \real^{n \times n}$ is called a {\em stable matrix},
%  if its spectral abscissa satisfies $\alpha(A) < 0$.
%\end{definition}
%%%%%%%%%%%%%%%%%%%%%%%%%%%%%%%%%%%%%%%%%%%%%%%%%%%%%%%%%%%%%%%%%%%%%%%%%%%%%

%%%%%%%%%%%%%%%%%%%%%%%%%%%%%%%%%%%%%%%%%%%%%%%%%%%%%%%%%%%%%%%%%%%%%%%%%%%%%
\begin{definition} \label{def:stable-pencil}
  A matrix pencil $(E,A)$ is called {\em stable}, if and only if
  each eigenvalue~$\lambda$ characterised by
  $\det (\lambda E - A) = 0$ has a strictly negative real part.
  % ${\rm Re}(\lambda)<0$.
\end{definition}
%%%%%%%%%%%%%%%%%%%%%%%%%%%%%%%%%%%%%%%%%%%%%%%%%%%%%%%%%%%%%%%%%%%%%%%%%%%%%

%%%%%%%%%%%%%%%%%%%%%%%%%%%%%%%%%%%%%%%%%%%%%%%%%%%%%%%%%%%%%%%%%%%%%%%%%%%%%
\begin{definition} \label{def:stable-system}
  The linear dynamical system~(\ref{linear-system}) is
  {\em asymptotically stable} if and only if its associated matrix pencil
  $(E,A)$ is stable.
\end{definition}
%%%%%%%%%%%%%%%%%%%%%%%%%%%%%%%%%%%%%%%%%%%%%%%%%%%%%%%%%%%%%%%%%%%%%%%%%%%%%

In the case of a non-singular mass matrix, asymptotic stability
of a system~(\ref{linear-system}) is equivalent to the property
$\alpha (E^{-1}A) < 0$ of the spectral abscissa in
Definition~\ref{def:spectral}.
Concerning Definition~\ref{def:stable-pencil},
a regular matrix pencil exhibits a finite set of eigenvalues.
Furthermore, Definition~\ref{def:stable-system} of asymptotic stability
can also be found in~\cite[p.~376]{braun}.
%The asymptotic stability guarantees the existence of the
%integral in~(\ref{h2-norm}).
The asymptotic stability guarantees the existence of the
transfer function~(\ref{transfer}) on the imaginary axis.
The $\htwo$-norm~(\ref{h2-norm}) is always finite for asymptotically
stable ODEs, whereas this norm may not exist in the case of (stable) DAEs.

If the matrix pencil $(E,A)$ has eigenvalues with non-positive
real part and a real part zero appears, then Lyapunov stability may still
be satisfied.
We consider this instance also as a loss of stability,
because the advantageous asymptotic stability is not valid any more.

%%%%%%%%%%%%%%%%%%%%%%%%%%%%%%%%%%%%%%%%%%%%%%%%%%%%%%%%%%%%%%%%%%%%%%%%%%%%%
\subsection{Projection-based model order reduction}
\label{sec:mor}
We assume that the linear dynamical system~(\ref{linear-system}) exhibits
a huge dimensionality~$n$.
Thus the involved matrices~$A$ and~$E$ are typically sparse.
The purpose of MOR is to decrease the complexity.
An alternative linear dynamical system
\begin{equation} \label{system-reduced}
  \begin{array}{rcl}
    \bar{E} \dot{\bar{x}}(t) & = & \bar{A} \bar{x}(t) + \bar{B} u(t) \\
    \bar{y}(t) & = & \bar{C} \bar{x}(t) \\
  \end{array}
\end{equation}
has to be constructed with state/inner variables
$\bar{x} : [0,\infty) \rightarrow \real^r$
and the matrices $\bar{A},\bar{E} \in \real^{r \times r}$,
$\bar{B} \in \real^{r \times n_{\rm in}}$, $\bar{C} \in \real^{n_{\rm out} \times r}$,
where the dimension~$r$ is much smaller than~$n$.
Initial values $\bar{x}(0)=\bar{x}_0$ are derived from
the initial values $x(0) = x_0$.
Nevertheless, the output of~(\ref{system-reduced}) should be a
good approximation to the output of~(\ref{linear-system}),
i.e., $\bar{y}(t) \approx y(t)$ for all relevant times.
The system~(\ref{system-reduced}) is called the reduced-order model (ROM)
of the full-order model (FOM) given by~(\ref{linear-system}).

The linear dynamical system~(\ref{system-reduced}) has its own
transfer function
$\bar{H} : \complex \backslash \bar{\Sigma} \rightarrow
\complex^{n_{\rm out} \times n_{\rm in}}$
of the form~(\ref{transfer}).
If both the original system~(\ref{linear-system}) and the
reduced system~(\ref{system-reduced}) are asymptotically stable,
then error bounds are available in the case of $x_0=0$ and $\bar{x}_0=0$.
% Assuming that the $\htwo$-norm~(\ref{h2-norm}) of $H - \bar{H}$ exists.
It holds that, see~\cite[p.~496]{benner-gugercin-willcox},
\begin{equation} \label{error-bound}
  \sup_{t \ge 0} \| y(t) - \bar{y}(t) \|_\infty \le
  \left\| H - \bar{H} \right\|_{\htwo} \| u \|_{\ltwo}
\end{equation} 
with the $\ltwo$-norm
\begin{equation} \label{l2-norm}
  \| u \|_{\ltwo} =
  \sqrt{ \int_0^{\infty} \| u(t) \|_2^2 \; {\rm d}t } \; , 
\end{equation}
the $\htwo$-norm~(\ref{h2-norm}),
the maximum (vector) norm $\| \cdot \|_\infty$ and the
Euclidean (vector) norm $\| \cdot \|_2$.

In projection-based MOR, see~\cite{antoulas},
 each approach yields
two projection matrices $V,W \in \real^{n \times r}$ of full rank.
We obtain the matrices of the ROM~(\ref{system-reduced}) by
\begin{equation} \label{projected-matrices}
  \bar{A} = W^\top A V , \quad \bar{B} = W^\top B , \quad
  \bar{C} = C V , \quad \bar{E} = W^\top E V .
\end{equation}
The orthogonality $V^\top V = I_r$ and sometimes the
biorthogonality $W^\top V = I_r$ are supposed
with the identity matrix $I_r \in \real^{r \times r}$.

Often the projection matrices result from the determination of
subspaces, i.e.,
$$ \mathcal{V} = {\rm span}(V) \subset \real^n
   \qquad \mbox{and} \qquad
   \mathcal{W} = {\rm span}(W) \subset \real^n . $$
On the one hand, the original state/inner variables are approximated
within the space~$\mathcal{V}$ by $x \approx V \bar{x}$.
On the other hand, the residual
\begin{equation} \label{residual}
  g(t) = E V \dot{\bar{x}}(t) - A V \bar{x}(t) - B u(t)
  \in \real^n
\end{equation}
is kept small by the requirement $g(t) \perp \mathcal{W}$
and thus $W^\top g(t) = 0$ for all~$t$.

%%%%%%%%%%%%%%%%%%%%%%%%%%%%%%%%%%%%%%%%%%%%%%%%%%%%%%%%%%%%%%%%%%%%%%%%%%%%%
\subsection{Galerkin-type methods}
\label{sec:galerkin}
A Galerkin-type projection~(\ref{projected-matrices}) is characterised
by the property~$W=V$.
Thus we have to determine just a suitable projection matrix~$V$.
Examples of Galerkin-type MOR methods are:
\begin{itemize}
\item
  one-sided Arnoldi method, see~\cite{freund},
\item
  proper orthogonal decomposition (POD), see~\cite[p.~277]{antoulas},
\item
  multi-parameter moment matching as in~\cite{li-etal},
\item
  iterative improvement for the case of many outputs
  as in~\cite{freitas},
\item
  and others.
\end{itemize}
Moment matching methods identify an approximation of the
transfer function~(\ref{transfer}) in the frequency domain.
The one-sided Arnoldi scheme represents a Galerkin-type
moment matching method. 
Alternatively, the POD technique employs information on a solution
for a particular input in the time domain.

We explain the one-sided Arnoldi method, because it is used for
the numerical experiments in Section~\ref{sec:examples}.
An expansion point $s_0 \in \complex \backslash \Sigma$ is chosen.
The matrix
\begin{equation} \label{matrix-moment}
    F = s_0 E - A \in \complex^{n \times n}
\end{equation}
is arranged, which includes the matrices of the
linear dynamical system~(\ref{linear-system}).
Let a single input ($n_{\rm in}=1$) be given without loss of generality.
We define the matrix $G = F^{-1} E$ and the vector $z = F^{-1} B$.
%Starting from $v_1 = - F^{-1} B$ and $\hat{v}_1 = \frac{1}{\| v_1 \|} v_1$,
%an orthogonalization procedure is applied.
The Krylov subspaces belonging to the matrix and the vector read as
\begin{equation} \label{krylov}
  \mathcal{K}_r (G,z) =
   {\rm span} \{ z,Gz,G^2z,\ldots,G^{r-1}z \} \subset \complex^n
\end{equation}
for $r \ge 1$.
The Arnoldi algorithm computes an orthonormal basis of the
subspace~(\ref{krylov}) by a specific orthogonalisation scheme.
The basis is collected in a matrix $\hat{V} \in \complex^{n \times r}$.
Hence it holds that $\hat{V}^\top \hat{V} = I_r$ by construction.
The projection matrix $V=\hat{V}$ becomes real-valued in the case of
$s_0 \in \real \backslash \Sigma$.
Otherwise, the projection matrix~$V$
is obtained from a reorthogonalisation of the matrix
$({\rm Re}(\hat{V}) , {\rm Im}(\hat{V})) \in \real^{n \times 2r}$.
This technique can be generalised straightforward to the case of
several expansion points and multiple inputs.

The computational effort of the one-sided Arnoldi method consists in
two parts.
Firstly, we compute an $LU$-decomposition of the
high-dimensional matrix~(\ref{matrix-moment}).
The sparsity of the matrices often allows for an efficient computation.
Each matrix-vector multiplication $Gz'$ with some $z' \in \complex^n$
requires a solution of a linear system with coefficient
matrix~(\ref{matrix-moment}).
Thus $r$ forward and backward substitutions are performed to
determine~(\ref{krylov}).
Secondly, basic linear algebra operations yield the orthonormal basis
similar to the Gram-Schmidt orthogonalisation.

In each MOR approach, the reduced system~(\ref{system-reduced}) is often
useless if it is not at least Lyapunov stable.
In particular, the error bound~(\ref{error-bound}) holds true only for
asymptotically stable systems.
Many moment matching techniques like Krylov subspace methods and POD
do not guarantee a stable ROM.

%%%%%%%%%%%%%%%%%%%%%%%%%%%%%%%%%%%%%%%%%%%%%%%%%%%%%%%%%%%%%%%%%%%%%%%%%%%%%
%%%                      Stability preservation                           %%%
%%%%%%%%%%%%%%%%%%%%%%%%%%%%%%%%%%%%%%%%%%%%%%%%%%%%%%%%%%%%%%%%%%%%%%%%%%%%%

%%%%%%%%%%%%%%%%%%%%%%%%%%%%%%%%%%%%%%%%%%%%%%%%%%%%%%%%%%%%%%%%%%%%%%%%%%%%%
\section{Stability preservation}
\label{sec:preservation}
We investigate linear dynamical systems~(\ref{linear-system})
consisting of ODEs in this section.

%%%%%%%%%%%%%%%%%%%%%%%%%%%%%%%%%%%%%%%%%%%%%%%%%%%%%%%%%%%%%%%%%%%%%%%%%%%%%
\subsection{Stability and Lyapunov equations}
We consider the Lyapunov inequality
\begin{equation}
  \label{lyapunov-ineq}
  A^\top M E + E^\top M A < 0 
\end{equation}
with $A,E$ from the dynamical system~(\ref{linear-system})
and a (non-unique) solution $M \in \real^{n \times n}$.
This problem consists in finding a symmetric positive definite
matrix~$M$ such that the left-hand side of~(\ref{lyapunov-ineq})
is negative definite.
We can solve the problem by choosing any symmetric positive definite matrix
$F \in \real^{n \times n}$.
It follows that the generalised Lyapunov equation,
cf.~\cite{penzl},
\begin{equation} \label{lyapunov}
  A^\top M E + E^\top M A + F = 0 
\end{equation}
yields a unique symmetric positive definite solution~$M$,
because the spectral abscissa from Definition~\ref{def:spectral}
exhibits $\alpha(E^{-1}A)<0$.
This solution~$M$ also satisfies the Lyapunov
inequality~(\ref{lyapunov-ineq}).

Direct methods of linear algebra compute the solution of~(\ref{lyapunov})
or its Cholesky factorisation, see~\cite{hammarling,penzl}.
However, direct methods are excluded in our context, since
their computational effort is $O(n^3)$.
Approximate methods are available like projection methods and
the alternating direction implicit (ADI) iteration,
see~\cite{kramer,li-white,wolf},
for example.
These methods often produce approximations $M \approx Z Z^\top$ with
a low-rank factor $Z \in \real^{n \times \ell}$ ($\ell \ll n$).
Thus the approximation becomes a singular matrix, which makes
the transformation dubious.
Ill-conditioned reduced matrices arise sometimes as shown in~\cite{pulch-arxiv}.

Using the solution of the Lyapunov equation~(\ref{lyapunov}),
we transform the ODEs~(\ref{linear-system}) into the equivalent system
\begin{equation} \label{ode-trafo}
  E^\top M E \dot{x}(t) = E^\top M A \dot{x} + E^\top M B u(t) .
\end{equation}
This transformation operates only in the image space and not in
the state space.
Stability preservation is given in a Galerkin-type MOR of
the equivalent system~(\ref{ode-trafo}).
We cite a theorem, whose proof can be found in~\cite{pulch-arxiv},
for example.

%%%%%%%%%%%%%%%%%%%%%%%%%%%%%%%%%%%%%%%%%%%%%%%%%%%%%%%%%%%%%%%%%%%%%%%%%%%%%
\begin{theorem} \label{thm:stable}
  Let the linear dynamical system~(\ref{linear-system}) with a non-singular
  mass matrix be asymptotically stable.
  If $M$ is a solution of the Lyapunov inequality~(\ref{lyapunov-ineq}),
  then each Galerkin-type projection-based MOR of
  the linear dynamical system~(\ref{ode-trafo}) 
  yields an asymptotically stable reduced system~(\ref{system-reduced}).
\end{theorem}
%%%%%%%%%%%%%%%%%%%%%%%%%%%%%%%%%%%%%%%%%%%%%%%%%%%%%%%%%%%%%%%%%%%%%%%%%%%%%

Let $V \in \real^{n \times r}$ with $V^\top V = I_r$ be any projection matrix
constructed for a reduction of the FOM~(\ref{linear-system}).
The Galerkin-type MOR with~$V$ is applied to the
transformed system~(\ref{ode-trafo}) now.
The matrices of the associated ROM~(\ref{system-reduced})
can be written in the form~(\ref{projected-matrices})
with the projection matrix
\begin{equation} \label{matrix-W}
  W = M E V .
\end{equation}
Consequently, this reduction represents a special case of a
(non-Galerkin-type) MOR~(\ref{projected-matrices}) for the
original system~(\ref{linear-system}).
It holds that $W^\top V \neq I_r$ in general.
Biorthogonality is achieved by the alternative projection matrix
\begin{equation} \label{matrix-tilde-W}
  W' = W (V^\top W)^{-1} .
\end{equation}
The ROMs obtained from~(\ref{matrix-W}) and~(\ref{matrix-tilde-W}) 
are equivalent and thus the stability properties coincide.
The matrix~$W$ in~(\ref{matrix-W}) is defined by matrix-matrix products.
On the one hand,
the evaluation $V'=EV$ is cheap, because~$E$ is typically sparse. 
On the other hand, the matrix~$M$ is dense.
We require an approximation of~$M$ such that the product $M V'$
is computable with relatively low effort.
The Galerkin-type MOR for the system~(\ref{linear-system}) and
the MOR for the system~(\ref{ode-trafo}) are not equivalent,
since the residual~(\ref{residual}) is not invariant in
the used transformations.

%%%%%%%%%%%%%%%%%%%%%%%%%%%%%%%%%%%%%%%%%%%%%%%%%%%%%%%%%%%%%%%%%%%%%%%%%%%%%
\subsection{Frequency domain integrals}
There are two analytical formulas for the solution of the
generalised Lyapunov equations~(\ref{lyapunov}),
see~\cite[p.~177]{antoulas}.
It holds that
\begin{equation} \label{lyap-time-domain}
  M = \int_0^\infty {\rm e}^{t (E^{-1} A)^\top}
  E^{-1} F E^{-\top} {\rm e}^{t (E^{-1}A)} \; {\rm d}t
\end{equation}
including the matrix exponential in the time domain.
The asymptotic stability of~(\ref{linear-system}) implies $\alpha(E^{-1}A)<0$
and thus this matrix-valued integral exists.
Alternatively, Parseval's theorem induces the matrix-valued integral
\begin{equation} \label{lyap-frequency-domain}
  M = \frac{1}{2\pi}
  \int_{-\infty}^\infty \big( {\rm i} \omega E^\top - A^\top \big)^{-1} F
  \big( -{\rm i} \omega E - A \big)^{-1} \; {\rm d}\omega
\end{equation}
in the frequency domain.
The asymptotic stability of~(\ref{linear-system}) yields
the invertibility of the involved matrices. % invertibility ! 
Although the integrand is complex-valued,
the integral~$M$ becomes real-valued.

Our idea is to apply the elementary choice~$F=I_n$ with the identity matrix
in the Lyapunov equation~(\ref{lyapunov}).
The identity matrix features the maximum rank~$n$.
Furthermore, there is no potential for a low-rank approximation of~$I_n$,
because all eigenvalues are identical to one.
A symmetry of the integrand allows for a restriction of the integration
to non-negative frequencies.
The frequency domain integral~(\ref{lyap-frequency-domain}) simplifies to
\begin{equation} \label{lyap-frequency-domain2}
  M = \frac{1}{\pi} \, {\rm Re} \left[ 
    \int_{0}^\infty S(\omega)^{-\htop} S(\omega)^{-1} \; {\rm d}\omega \right]
\end{equation}
using the abbreviation
\begin{equation} \label{matrix-S}
  S(\omega) = -{\rm i} \omega E - A \in \complex^{n \times n}
  \qquad \mbox{for}\;\; \omega \in \real . 
\end{equation}
In~(\ref{lyap-time-domain}), the matrix exponential yield dense matrices.
In~(\ref{lyap-frequency-domain}), the matrices $s E - A$
for $s \in \complex$ are often sparse,
whereas the inverse matrices are always dense.
Hence we never compute the inverse matrices explicitly.
Nevertheless, sophisticated algorithms often produce a sparse
$LU$-decomposition of a matrix $s E - A$.

In an MOR with projection matrix~(\ref{matrix-W}),
we require just the matrix-matrix product $M V'$
with the constant matrix $V' = EV$.
It follows that
\begin{equation} \label{M-times-V}
  W = M V' =  \frac{1}{\pi} \, {\rm Re} \left[ 
    \int_{0}^\infty S(\omega)^{-\htop} S(\omega)^{-1} V' \; {\rm d}\omega
    \right] ,
\end{equation}
which represents a matrix-valued integral of size $n \times r$
in the frequency domain.

We do not use the formulation~(\ref{lyap-time-domain}) in the time domain,
because there is no suitable method to calculate the
matrix exponential, see~\cite{moler-vanloan}.
Appropriate iterative techniques to compute a matrix-vector product
with the matrix exponential do exist, see~\cite{almohy-higham}.
However, rough approximations would confuse the error estimation of
an adaptive quadrature method.

%%%%%%%%%%%%%%%%%%%%%%%%%%%%%%%%%%%%%%%%%%%%%%%%%%%%%%%%%%%%%%%%%%%%%%%%%%%%%
\subsection{Error condition for Lyapunov inequality}

We show the following general result to characterise the influence of
errors in the context of the Lyapunov inequality~(\ref{lyapunov-ineq}).

%%%%%%%%%%%%%%%%%%%%%%%%%%%%%%%%%%%%%%%%%%%%%%%%%%%%%%%%%%%%%%%%%%%%%%%%%%%%%
\begin{theorem} \label{thm:condition}
  Let $M \in \real^{n \times n}$ be the solution of the
  Lyapunov equations~(\ref{lyapunov})
  with $F=I_n$ and $\widetilde{M} \in \real^{n \times n}$
  be any symmetric matrix.
  If it holds that
  \begin{equation} \label{error-tolerance-abs}
   \| \widetilde{M} - M \| <  \displaystyle
   \frac{1}{ \| A^\top \| \cdot \| E \| +  \| A \| \cdot \| E^\top \|}
  \end{equation}
  %\begin{eqnarray} \label{error-tolerance-abs}
  %    \| \widetilde{M} - M \| & < & \displaystyle
  %    \frac{1}{ \| A^\top \| \cdot \| E \| +  \| A \| \cdot \| E^\top \|},
  %    \quad \mbox{or} \\
  %    \displaystyle \label{error-tolerance-rel}
  %     \frac{\| \widetilde{M} - M \|}{\| M \|} & < & 1  
  %\end{eqnarray}
  in some subordinate matrix norm~$\| \cdot \|$,
  then $\widetilde{M}$ is positive definite and satisfies
  the Lyapunov inequality~(\ref{lyapunov-ineq}).
\end{theorem}
%%%%%%%%%%%%%%%%%%%%%%%%%%%%%%%%%%%%%%%%%%%%%%%%%%%%%%%%%%%%%%%%%%%%%%%%%%%%%

\underline{Proof:}

It holds that
$$ \begin{array}{rcl}
A^\top M E + E^\top M A & = & - I_n \\[0.5ex] 
A^\top \widetilde{M} E + E^\top \widetilde{M} A & = & - G \\
\end{array} $$
with a symmetric matrix~$G$.
Subtraction yields
$$ A^\top (\widetilde{M}-M) E + E^\top (\widetilde{M}-M) A = I_n - G . $$
Let $\eta = \| A^\top \| \cdot \| E \| +  \| A \| \cdot \| E^\top \|$.
%We obtain
%$$ 1 = \| I_n \| = \| A^\top M E + E^\top M A \| \le 
%\left( \| A^\top \| \cdot \| E \| +  \| A \| \cdot \| E^\top \| \right)
%\| M \|
%=  \eta \| M \| $$
%and thus $1 \le \eta \| M \|$.
We estimate
$$ \| I_n - G \| =
\| A^\top (\widetilde{M}-M) E + E^\top (\widetilde{M}-M) A \| \le
\eta \| \widetilde{M}-M \| . $$
Now the condition~(\ref{error-tolerance-abs}) 
is sufficient for $\eta \| \widetilde{M}-M \| < 1$ and thus
$\| I_n - G \| < 1$. 

Let $\lambda_j$ and $\mu_j$ for $j=1,\ldots,n$ be the eigenvalues of~$G$
and $I_n - G$, respectively.
It follows that
$$ | 1 - \lambda_j | = | \mu_j | 
\le \| I_n - G \| < 1 \quad \mbox{for}\;\; j=1,\ldots,n . $$
We obtain $0< \lambda_j < 2$ for all $j=1,\ldots,n$.
Consequently, the matrix~$G$ is positive definite
and the matrix~$-G$ is negative definite.
Since $\widetilde{M}$ represents the solution of a
Lyapunov equation~(\ref{lyapunov}) including
the symmetric positive definite matrix~$F=G$, $\widetilde{M}$~inherits the
positive definiteness.
\hfill $\Box$

\medskip

If we employ the spectral (matrix) norm, then the above constant~$\eta$
simplifies to $\eta = 2 \| A \|_2 \| E \|_2$.
However, the evaluation of the spectral norm takes more effort
in comparison to the norms $\|\cdot\|_1,\|\cdot\|_\infty$.
%Furthermore, a lower bound on the solution reads as
%$\frac{1}{\eta} \le \| M \|$ in any subordinate norm.

An obvious question is if a sufficient condition can be derived
for the relative error $\frac{\| \widetilde{M} - M \|}{\| M \|}$.
However, we require an upper bound on $\| M \|$ in this case,
which becomes more involved.
In~\cite[p.~100]{stykel-diss}, the analysis yields the estimate
$$ \| M \|_{\rm F} = \| \mathcal{L}^{-1} (I_n) \|_{\rm F} \le
\| \mathcal{L}^{-1} \|_{\rm F} \| I_n \|_{\rm F} = \sqrt{n}
\left( \inf_{\| X \|_{\rm F} = 1} \| A^T X E + E^\top X A \|_{\rm F} \right)^{-1} $$
including the inverse of the Lyapunov operator $\mathcal{L}$
in the Frobenius norm.
Yet this upper bound cannot be simplified in the case of general
matrices~$A$ and~$E$.

We obtain a necessary condition on the relative error with respect
to the requirement~(\ref{error-tolerance-abs}).

%%%%%%%%%%%%%%%%%%%%%%%%%%%%%%%%%%%%%%%%%%%%%%%%%%%%%%%%%%%%%%%%%%%%%%%%%%%%%
\begin{lemma} \label{lemma:relative-error}
  If the solution~$M$ of the Lyapunov equations~(\ref{lyapunov})
  with $F=I_n$ and a symmetric matrix~$\widetilde{M}$ satisfy the
  condition~(\ref{error-tolerance-abs}),
  then the relative error exhibits the bound
  \begin{equation} \label{delta-M-rel}
    \frac{\| \widetilde{M} - M \|}{\| M \|} < 1
  \end{equation}
  in the used matrix norm.
\end{lemma}
%%%%%%%%%%%%%%%%%%%%%%%%%%%%%%%%%%%%%%%%%%%%%%%%%%%%%%%%%%%%%%%%%%%%%%%%%%%%%

\underline{Proof:}

The Lyapunov equation~(\ref{lyapunov}) with $F=I_n$ yields
$$ 1 = \| I_n \| \le \| A^\top M E + E^\top M A \| \le \eta \| M \| $$
with the constant~$\eta$ in the proof of Theorem~\ref{thm:condition}.
It follows that $\frac{1}{\eta} \le \| M \|$. 
We obtain 
$$ \frac{\| \widetilde{M} - M \|}{\| M \|} \le
\frac{\| \widetilde{M} - M \|}{\frac{1}{\eta}} =
\eta \| \widetilde{M} - M \| < 1 $$
using the property~(\ref{error-tolerance-abs}).
\hfill $\Box$

\medskip

Lemma~\ref{lemma:relative-error} motivates that the
condition~(\ref{error-tolerance-abs}) is not strong,
because the induced relative error~(\ref{delta-M-rel})
may become up to 100\%.

%%%%%%%%%%%%%%%%%%%%%%%%%%%%%%%%%%%%%%%%%%%%%%%%%%%%%%%%%%%%%%%%%%%%%%%%%%%%%
\subsection{Quadrature methods}
\label{sec:quadrature}
Phillips and Silveira~\cite{phillips-silveira} investigated
the integrals~(\ref{lyap-frequency-domain}) including
$F = G G^\top$ with a low-rank factor $G \in \real^{n \times \ell}$ ($\ell \ll n$).
Therein, a quadrature method with $K$~nodes and positive weights yields
an approximation $M \approx Z Z^{\htop}$ with a factor
$Z \in \complex^{n \times (\ell K)}$.
Hence both the rank~$\ell$ and the number of nodes~$K$ has to be small.
This requirement disappears in our approach using $F=I_n$.

We also apply a quadrature rule to our frequency domain integrals.
For theoretical investigations, we define an
approximation of~(\ref{lyap-frequency-domain2}) by
\begin{equation} \label{approx-M}
  \widetilde{M} =  \frac{1}{\pi} \, {\rm Re} \left[
    \sum_{k=1}^K \gamma_k S(\omega_k)^{-\htop} S(\omega_k)^{-1}
    \right] =
    \frac{1}{\pi} 
    \sum_{k=1}^K \gamma_k {\rm Re} \left[ S(\omega_k)^{-\htop} S(\omega_k)^{-1}
    \right]
\end{equation}
with nodes~$\omega_k \ge 0$ and weights~$\gamma_k > 0$ for $k=1,\ldots,K$.
The quadrature introduces an error.
Nevertheless, the approximations~(\ref{approx-M}) always own the
following desired properties independent of the magnitude of the error.

%%%%%%%%%%%%%%%%%%%%%%%%%%%%%%%%%%%%%%%%%%%%%%%%%%%%%%%%%%%%%%%%%%%%%%%%%%%%%
\begin{lemma} \label{lemma:definiteness}
  If the quadrature rule involves positive weights only, then the
  approximation~$\widetilde{M}$ in~(\ref{approx-M}) is always symmetric
  and positive definite.
\end{lemma}
%%%%%%%%%%%%%%%%%%%%%%%%%%%%%%%%%%%%%%%%%%%%%%%%%%%%%%%%%%%%%%%%%%%%%%%%%%%%%

\underline{Proof:} \nopagebreak

We define $S(\omega_k)^{-1} = X_k + {\rm i} Y_k$ with
real-valued matrices $X_k,Y_k$.
It follows that
$$ {\rm Re}  \left[ S(\omega_k)^{-\htop} S(\omega_k)^{-1} \right] =
X_k^{\top} X_k + Y_k^{\top} Y_k . $$
Now the symmetry of~$\widetilde{M}$ is obvious.
We show the definiteness.
Let $z \in \real^n \backslash \{ 0 \}$.
We obtain
$$ z^\top \widetilde{M} z =  \frac{1}{\pi} \sum_{k=1}^K \gamma_k
\left( z^{\top} X_k^{\top} X_k z + z^{\top} Y_k^{\top} Y_k z \right) =
 \frac{1}{\pi} \sum_{k=1}^K \gamma_k \left( \| X_k z \|_2^2 + \| Y_k z \|_2^2
\right) \ge 0 . $$
It holds that $X_k z \neq 0$ or $Y_k z \neq 0$ for each~$k$,
because otherwise $S(\omega_k)^{-1} z = 0$ would cause a contradiction
to the non-singularity of the matrix $S(\omega_k)$.
It follows that the above sum is strictly positive.
\hfill $\Box$

\medskip

%Since the approximation~(\ref{approx-M}) is symmetric,
%Theorem~\ref{thm:condition} yields a criterion to guarantee the
%preservation of stability.
%However, this error condition also cannot be checked for
%a quadrature rule in practise.

The associated approximation of~(\ref{M-times-V}) reads as
\begin{equation} \label{approx-W}
  \widetilde{W} =
  \widetilde{M} V' =
  %\frac{1}{\pi} \, {\rm Re} \left[
  %\sum_{k=1}^K \gamma_k S(\omega_k)^{-1} S(\omega_k)^{-\htop} V' \right] 
  \frac{1}{\pi} \sum_{k=1}^K \gamma_k \,
       {\rm Re} \left[ S(\omega_k)^{-\htop} S(\omega_k)^{-1} V' \right] .
\end{equation}
It turned out that this approach is similar to a quadrature technique
given in~\cite{benner-schneider}, where an integral of the
form~(\ref{lyap-frequency-domain}) yields the Gramian of a
linear dynamical system with many outputs.

Lemma~\ref{lemma:definiteness} shows that the matrix~$\widetilde{M}$
is always symmetric and positive definite in~(\ref{approx-W}).
Consequently, the underlying transformation is non-singular,
which represents a crucial advantage in comparison to methods for
Lyapunov equations~(\ref{lyapunov}) producing low-rank approximations.
Theorem~\ref{thm:condition} implies that an approximation
satisfying~(\ref{error-tolerance-abs}) guarantees
a stability preservation in an MOR.
However, the approximation errors cannot be checked
in practise, because the exact solution~$M$ is unknown.
We compute an approximation~(\ref{approx-W}), where an adaptive
quadrature yields errors below predetermined tolerances.
Yet there is no direct connection to the errors in the 
counterpart~(\ref{approx-M}).

%%%%%%%%%%%%%%%%%%%%%%%%%%%%%%%%%%%%%%%%%%%%%%%%%%%%%%%%%%%%%%%%%%%%%%%%%%%%%
\subsection{Numerical solution of linear systems}
Our aim is to evaluate the approximation~(\ref{approx-W})
without computations of matrices of size $n \times n$.
Thus we solve complex-valued linear systems 
\begin{equation} \label{linear-system-alg}
  S(\omega_k) S(\omega_k)^{\htop} X_k = V'
\end{equation}
for each~$k$
with the matrices~(\ref{matrix-S}), a predetermined matrix
$V' \in \real^{n \times r}$ and the unknowns $X_k \in \complex^{n \times r}$.
We consider only direct approaches of numerical linear algebra.
Iterative methods introduce an additional error,
which restricts the accuracy of high-order quadrature rules.
Our aim is to use as less nodes as possible.

There are two possibilities to solve a
linear system~(\ref{linear-system-alg}) directly now:
\begin{itemize}
\item[(i)]
  The matrix-matrix product
  $\hat{S}_k = S(\omega_k) S(\omega_k)^{\htop}$ is computed.
  An algorithm for sparse matrices generates a
  Cholesky-decomposition $\hat{S}_k = \hat{L}_k \hat{L}_k^{\htop}$.
  We determine the solution
  $X_k = \hat{L}_k^{-\htop} \hat{L}_k^{-1} V'$
  by forward and backward substitutions for multiple right-hand sides.
\item[(ii)]
  An $LU$-decomposition including pivoting with row as well as
  column reordering is applied to the matrices~(\ref{matrix-S}), i.e.,
  \begin{equation} \label{lu-decomp}
    P_k S(\omega_k) Q_k = L_k U_k
  \end{equation}
  with orthogonal permutation matrices~$P_k,Q_k$.
  % $P_k,Q_k \in \real^{n \times n}$
  It follows that
  $$ Q_k^{\top} S(\omega_k)^{\htop} P_k^{\top} = U_k^{\htop} L_k^{\htop} $$
  represents an $LU$-decomposition of $S(\omega_k)^{\htop}$.
  We obtain
  $$ X_k = P_k^\top L_k^{-\htop} U_k^{-\htop} U_k^{-1} L_k^{-1} P_k V' $$
  using permutations, forward and backward substitutions
  for multiple right-hand sides.
\end{itemize}
We do not apply the approach~(i) in our numerical computations for
two reasons:
(i)~The matrix $\hat{S}_k$ is less sparse than a matrix~(\ref{matrix-S}).
Thus a Cholesky factorisation of $\hat{S}_k$ may not be (significantly)
faster than an $LU$-decomposition of $S(\omega_k)$.
(ii)~The condition number increases considerably due
to ${\rm cond}(\hat{S}_k) = {\rm cond}(S(\omega_k))^2$
with respect to the spectral norm.

Hence we perform the $LU$-decompositions~(\ref{lu-decomp}).
Efficient numerical methods are available like UMFPACK~\cite{davis}.
Therein, pivoting and permutation strategies keep the factorisations
as sparse as possible,
while still numerical stability is achieved.

%%%%%%%%%%%%%%%%%%%%%%%%%%%%%%%%%%%%%%%%%%%%%%%%%%%%%%%%%%%%%%%%%%%%%%%%%%%%%
\subsection{Integrals on finite intervals}
\label{sec:finite-interval}
Concerning the integrals~(\ref{M-times-V}),
we can straightforward transform the infinite frequency domain $[0,\infty)$
to the finite interval $[0,1)$.
The substitution $\omega = \frac{\xi}{1-\xi}$
or, equivalently, $\xi = \frac{\omega}{1+\omega}$ yields
\begin{equation} \label{integral-trafo}
  W =  \frac{1}{\pi} \, {\rm Re} \left[ 
    \int_{0}^1 S\left(\frac{\xi}{1-\xi}\right)^{-\htop}
    S\left(\frac{\xi}{1-\xi}\right) V' \frac{1}{(1-\xi^2)} \; {\rm d}\xi
    \right] .
\end{equation}
An advantage is that the evaluation of the integrand at $\xi = 1$
exists in the limit
\begin{equation} \label{limit-xi}
  \lim_{\xi \rightarrow 1} \; S\left(\frac{\xi}{1-\xi}\right)^{-\htop}
  S\left(\frac{\xi}{1-\xi}\right) \frac{1}{(1-\xi^2)} = E^{-\top} E^{-1}
\end{equation}
provided that the mass matrix is non-singular.
Now any (open or closed) quadrature rule for finite intervals generates
an approximation to the integral~(\ref{integral-trafo}).
Numerical tests indicate that the Gauss-Legendre quadrature
is superior for computing the integral~(\ref{integral-trafo})
in comparison to other common schemes.
The reason is that the integrand is analytic in the open interval~$(0,1)$.

An adaptive Gauss-Kronrod quadrature rule was considered for integrals
of this type in~\cite{benner-schneider}.
However, as mentioned in Section~\ref{sec:quadrature},
the quadrature is not required to be sufficiently accurate for
the integrals~(\ref{integral-trafo}) but the inherent
integrals~(\ref{lyap-frequency-domain2}).
An alternative is to use an adaptive refinement of nested
quadrature rules until a desired set of ROMs becomes asymptotically stable.
Note that it is cheap to check the stability for low-dimensional systems.

An elementary adaptive method of this type can be based on
the midpoint rule.
The nodes read as $\xi_k = \frac{h}{2} + (k-1)h \in (0,1)$ for $k=1,\ldots,K$
with step size $h = \frac{1}{K}$.
The iteration $K_i = 2^{i-1}$ for $i=1,2,3,\ldots$ induces a sequence of
nested grids, where the evaluations of the integrand
in~(\ref{integral-trafo}) can be reused.
The iteration is terminated if all considered ROMs become stable.

%%%%%%%%%%%%%%%%%%%%%%%%%%%%%%%%%%%%%%%%%%%%%%%%%%%%%%%%%%%%%%%%%%%%%%%%%%%%%
%%%               Differential-Algebraic Equations                        %%%
%%%%%%%%%%%%%%%%%%%%%%%%%%%%%%%%%%%%%%%%%%%%%%%%%%%%%%%%%%%%%%%%%%%%%%%%%%%%%

\section{Application to differential-algebraic equations}
\label{sec:daes}
The technique of Section~\ref{sec:preservation} cannot be applied directly
to systems of DAEs~(\ref{linear-system}).
A singular mass matrix implies 
$z^\top (A^\top M E + E^\top M A) z=0$ for $z \in {\rm ker}(E)$ and any~$M$.
Hence the Lyapunov equation~(\ref{lyapunov}) is not fulfilled
for each definite matrix~$F$.
The associated integral~(\ref{lyap-frequency-domain}) does not exist,
even though the integrand is always well-defined in the case of
a stable matrix pencil with respect to Definition~\ref{def:stable-pencil}.
Likewise, the limit~(\ref{limit-xi}) does not exist.

%%%%%%%%%%%%%%%%%%%%%%%%%%%%%%%%%%%%%%%%%%%%%%%%%%%%%%%%%%%%%%%%%%%%%%%%%%%%%
\subsection{Kronecker normal form}
A linear dynamical system~(\ref{linear-system}) with a singular
mass matrix can be transformed into the Kronecker normal form,
see~\cite[p.~452]{hairer2}.
There are non-singular matrices $T_{\rm l},T_{\rm r} \in \real^{n \times n}$
such that
\begin{equation} \label{kronecker}
  T_{\rm l} A T_{\rm r} = \begin{pmatrix} A' & 0 \\ 0 & I_{n_2} \\ \end{pmatrix}
  \qquad \mbox{and} \qquad
  T_{\rm l} E T_{\rm r} = \begin{pmatrix} I_{n_1} & 0 \\ 0 & N \\ \end{pmatrix}
\end{equation}
with a matrix $A' \in \real^{n_1 \times n_1}$ and a nilpotent strictly upper
triangular matrix $N \in \real^{n_2 \times n_2}$ ($n=n_1+n_2$).
The system~(\ref{linear-system}) splits into a slow and a fast subsystem
\begin{equation} \label{dae-semi-expl}
  \begin{array}{rcl}
    \dot{z}_1(t) & = & A' z_1(t) \\
    N \dot{z}_2(t) & = & z_2(t) \\
  \end{array}
\end{equation}
with $z_1(t) \in \real^{n_1}$ and $z_2(t) \in \real^{n_2}$.
The input terms are omitted in~(\ref{dae-semi-expl}), because
they do not influence the stability properties of the systems.
The smallest integer~$\nu \ge 0$ such that $N^{\nu-1} \neq 0$ and $N^{\nu}=0$
is called the nilpotency index of the DAE system.

Unfortunately, there is no efficient numerical method to compute the
decomposition~(\ref{kronecker}) of the matrices in the
linear dynamical system~(\ref{linear-system}).
Thus we must design techniques, which are feasible without an
explicit knowledge of the Kronecker normal form.

%%%%%%%%%%%%%%%%%%%%%%%%%%%%%%%%%%%%%%%%%%%%%%%%%%%%%%%%%%%%%%%%%%%%%%%%%%%%%
\subsection{Regularisation}
\label{sec:regularisation}
In~\cite{mohaghegh}, a system of DAEs was regularised straightforward
under the assumption of semi-explicit systems.
Alternatively, we apply an approach for general descriptor systems
from~\cite{mueller}, which goes back to~\cite{wang-etal}.
The matrix pencil $(E,A)$ is modified into $(\hat{E},\hat{A})$ by
\begin{equation} \label{matrices-regularised}
  \hat{E} = E - \alpha A
  \qquad \mbox{and} \qquad
  \hat{A} = A + \beta E
\end{equation}
with parameters $\alpha,\beta > 0$.
The matrix $\hat{E}$ is non-singular for all $\alpha > 0$,
because otherwise some $\lambda = \frac{1}{\alpha} > 0$ would be
an eigenvalue of the stable matrix pencil $(E,A)$.
The theorem below follows from the results in~\cite{mueller}.

%%%%%%%%%%%%%%%%%%%%%%%%%%%%%%%%%%%%%%%%%%%%%%%%%%%%%%%%%%%%%%%%%%%%%%%%%%%%%
\begin{theorem} \label{thm:regularisation}
  Let the matrix~$E$ be singular and the matrix pencil $(E,A)$ be regular. 
  The perturbed matrices~(\ref{matrices-regularised})
  with $\alpha = \beta^2$ and
  \begin{equation} \label{bound-beta}
    0 < \beta < \frac{1}{\rho(A')^2}
  \end{equation}
  using the spectral radius $\rho(A')$ of the matrix
  in the Kronecker normal form~(\ref{kronecker})
  yield an asymptotically stable system of ODEs~(\ref{linear-system}).
\end{theorem}
%%%%%%%%%%%%%%%%%%%%%%%%%%%%%%%%%%%%%%%%%%%%%%%%%%%%%%%%%%%%%%%%%%%%%%%%%%%%%

The upper bound~(\ref{bound-beta}) is unknown in general,
because the Kronecker normal form~(\ref{kronecker})
is not available in practise.
Nevertheless, Theorem~\ref{thm:regularisation} tells us that a
regularisation to a stable system is feasible for all
sufficiently small~$\beta>0$.
Furthermore, the sparsity pattern of $s E - A$ is identical to
the sparsity pattern of $s \hat{E} - \hat{A}$.
Thus the computational effort for a solution of linear systems
does not change significantly.

Now we employ MOR methods to the regularised system~(\ref{linear-system})
including the matrices~(\ref{matrices-regularised}),
where the stability-preserving technique is applicable
from Section~\ref{sec:preservation}.

%%%%%%%%%%%%%%%%%%%%%%%%%%%%%%%%%%%%%%%%%%%%%%%%%%%%%%%%%%%%%%%%%%%%%%%%%%%%%
\subsection{Error estimates}
We reduce the regularised system with matrices~(\ref{matrices-regularised})
instead of the original descriptor system~(\ref{linear-system}).
This approach is appropriate, if error bounds can be provided for
the transfer functions.
It holds that
\begin{equation} \label{dae-total-error}
  \left\| H_{\rm DAE} - H_{\rm ROM} \right\| \le
  \left\| H_{\rm DAE} - H_{\rm ODE} \right\| +
  \left\| H_{\rm ODE} - H_{\rm ROM} \right\|
\end{equation}
in each norm, where $H_{\rm ODE}$ and $H_{\rm ROM}$ are the transfer functions
of the regularised system and its ROM, respectively.
The difference $H_{\rm ODE} - H_{\rm ROM}$ depends on the quality of the MOR.
We discuss the difference $H_{\rm DAE} - H_{\rm ODE}$ in the following.

A linear dynamical system (or its transfer function)
is called strictly proper,
if and only if the transfer function~$H$ satisfies
$$ \lim_{s \rightarrow \infty} H(s) = 0 . $$
A system of ODEs always exhibits a strictly proper transfer function.
The transfer function of a general system of DAEs reads as,
see~\cite{benner-stykel},
\begin{equation} \label{transfer-dae}
  H_{\rm DAE}(s) = H_{\rm SP}(s) + P(s)
\end{equation}
with a strictly proper part
$H_{\rm SP} : \complex \backslash \Sigma \rightarrow
\complex^{n_{\rm out} \times n_{\rm in}}$
and a polynomial part
$P : \complex \rightarrow \complex^{n_{\rm out} \times n_{\rm in}}$.
The polynomial part either vanishes or represents a non-zero matrix-valued
polynomial of degree at most~$\nu$ with the index~$\nu$ of the system.
These properties also depend on the definition of inputs and outputs
in each system.

If the polynomial part vanishes, then the $\htwo$-norm~(\ref{h2-norm})
as well as the $\hinf$-norm of the transfer function~(\ref{transfer-dae})
exist independent of the index.
The $\hinf$-norm is always finite in the case of index-one systems. 
In~\cite{guenther}, for example,
the electric circuit of the Miller integrator is modelled by a
linear DAE system of index $\nu=2$, where the polynomial part becomes zero
for the transfer function relating the input voltage to the output voltage.

We provide an error bound for the regularisation
on compact frequency intervals,
which is relevant in this context.

%%%%%%%%%%%%%%%%%%%%%%%%%%%%%%%%%%%%%%%%%%%%%%%%%%%%%%%%%%%%%%%%%%%%%%%%%%%%%
\begin{lemma} \label{lemma:error}
  Assume that the DAE is asymptotically stable and the ODE is given
  by~(\ref{matrices-regularised}) with $\alpha = \beta^2$.
  For each $\omega' > 0$ there are constants
  $K_{\omega'},L_{\omega'} > 0$ such that
  \begin{equation} \label{error-compact}
    \left\| H_{\rm DAE}({\rm i}\omega) - H_{\rm ODE}({\rm i}\omega) \right\|_2
    < K_{\omega'} \beta
  \end{equation}
  uniformly for all frequencies $\omega \in [-\omega',\omega']$
  provided that $\beta<L_{\omega'}$,
  where the spectral (matrix) norm $\| \cdot \|_2$ is used.
\end{lemma}
%%%%%%%%%%%%%%%%%%%%%%%%%%%%%%%%%%%%%%%%%%%%%%%%%%%%%%%%%%%%%%%%%%%%%%%%%%%%%

\underline{Proof:}

Let $I = [-\omega',\omega']$, $s={\rm i}\omega$ and
$\| \cdot \| = \| \cdot \|_2$ in this proof.
We assume a bound
$L_{\omega'} \le \min \{ \rho(A')^{-2} , 1 \}$ on~$\beta$
due to~(\ref{bound-beta}).
Theorem~\ref{thm:regularisation} implies that the ODE is
asymptotically stable.
Thus the transfer functions exist and are continuous on the imaginary axis.
The spectral norm is submultiplicative for matrices of any size.
Thus we estimate
$$ \left\| H_{\rm DAE}(s) - H_{\rm ODE}(s) \right\| \le
\| B \| \cdot \| C \| \cdot
\left\| (sE-A)^{-1} - (s\hat{E}-\hat{A})^{-1} \right\| . $$
We use the abbreviations
$G(s) = s E - A$ and $\hat{G}(s) = s \hat{E} - \hat{A}$.
A general estimate on the difference between inverse matrices is available
in a subordinate matrix norm.
It follows that
$$ \left\| G(s)^{-1} - \hat{G}(s)^{-1} \right\| \le
\frac{\| G(s)^{-1} \|^2 \| G(s) - \hat{G}(s) \|}{
  1- \| G(s)^{-1} \| \cdot \| G(s) - \hat{G}(s) \|} $$
provided that $\| G^{-1}(s) \| \cdot \| G(s) - \hat{G}(s) \| < 1$. 
On the one hand, the definition~(\ref{matrices-regularised})
with $\alpha = \beta^2$ and $\beta \le 1$ implies
$$  \left\| G(s) - \hat{G}(s) \right\| =
\| s \alpha A - \beta E \| \le
\beta^2 |s| \cdot \| A \| + \beta \| E \| \le
\Gamma \beta $$
for all $\omega \in I$ with the constant
$\Gamma = \omega' \| A \| + \| E \| > 0$.
On the other hand, we require the constant
$$ \Theta = \max_{\omega \in I}
\left\| \left( {\rm i}\omega E - A \right)^{-1} \right\| > 0 $$
such that $\| G(s)^{-1} \| \le \Theta$ for all $\omega \in I$.
We obtain
$$ \left\| G(s)^{-1} - \hat{G}(s)^{-1} \right\| \le
2 \Theta^2 \Gamma \beta $$
for all $\omega \in I$ provided that
$\beta < \frac{1}{2 \Theta \Gamma}$. 
Consequently, the constants read as
$K_{\omega'} := 2 \Theta^2 \Gamma \| B \| \| C \|$ and
$L_{\omega'} = \min \{ \frac{1}{2\Theta\Gamma} , \frac{1}{\rho(A')^2} , 1\}$.
\hfill $\Box$

\medskip

Lemma~\ref{lemma:error} demonstrates that error of the regularisation
is low on a compact frequency domain
for sufficiently small parameters $\alpha,\beta$.
An error estimate of the type~(\ref{error-compact}) cannot be derived
uniformly for all $\omega \in \real$,
because the integral~(\ref{lyap-frequency-domain}) does not exist
in the limit $\beta \rightarrow 0$.
Thus high frequencies represent the critical part.
If a system of DAEs~(\ref{linear-system}) has a strictly proper
transfer function, then this problem becomes obsolete.

%%%%%%%%%%%%%%%%%%%%%%%%%%%%%%%%%%%%%%%%%%%%%%%%%%%%%%%%%%%%%%%%%%%%%%%%%%%%%
\begin{theorem} \label{thm:error}
  Let the linear dynamical system~(\ref{linear-system}) be an
  asymptotically stable DAE with a strictly proper transfer function.
  A system of ODEs is given by the perturbed
  matrices~(\ref{matrices-regularised}) with $\alpha = \beta^2$.
  For each $\varepsilon > 0$ there is a constant $L > 0$ such that
  the transfer functions satisfy
  \begin{equation} \label{error-global}
    \left\| H_{\rm DAE} - H_{\rm ODE} \right\|_{\htwo} < \varepsilon
  \end{equation}
  for all parameters $\beta$ with $0<\beta<L$.
\end{theorem}
%%%%%%%%%%%%%%%%%%%%%%%%%%%%%%%%%%%%%%%%%%%%%%%%%%%%%%%%%%%%%%%%%%%%%%%%%%%%%

\underline{Proof:}

The restriction $\beta \le \beta' < \min \{ \rho(A')^{-2},1 \}$,
see~(\ref{bound-beta}), for some $\beta'>0$ (close to the upper bound)
and $\alpha = \beta^2$ ensures the existence of each $\htwo$-norm.
The $\htwo$-norm of the difference reads as
$$ \left\| H_{\rm DAE} - H_{\rm ODE} \right\|_{\htwo}^2 =
\frac{1}{2\pi}
\int_{-\infty}^{\infty}
\left\| H_{\rm DAE}({\rm i}\omega) - H_{\rm ODE}({\rm i}\omega) \right\|_{\rm F}^2
\; {\rm d}\omega $$
including the Frobenius (matrix) norm.
Let $\Delta H = H_{\rm DAE} - H_{\rm ODE}$.
The assumptions imply that each component of
$\Delta H$ is a rational function of~$\omega$,
where the degree of the numerator is less than the degree of the
denominator.
The coefficients of the rational functions depend continuously on
the parameters~$\alpha,\beta \ge 0$.

We discuss the part for high frequencies, i.e.,
$$ \int_{\omega'}^{\infty} \left\| \Delta H ({\rm i}\omega) \right\|_{\rm F}^2
\; {\rm d}\omega =
\sum_{i=1}^{n_{\rm in}} \sum_{j=1}^{n_{\rm out}} \int_{\omega'}^{\infty}
\left| \Delta H_{ij} ({\rm i}\omega) \right|^2 \;{\rm d}\omega $$
with $\omega' \gg 1$.
It follows that
$$ \lim_{\omega' \rightarrow \infty} \int_{\omega'}^{\infty}
 \left\| \Delta H ({\rm i}\omega) \right\|_{\rm F}^2
 \; {\rm d}\omega = 0 $$
for each $\beta \in \mathcal{B} = [0,\beta']$ and $\alpha = \beta^2$,
where the convergence is monotone from above
(for $\omega' \rightarrow \infty$).
The parameter interval $\mathcal{B}$ is compact.
Dini's theorem yields
$$  \lim_{\omega' \rightarrow \infty} \max_{\beta \in \mathcal{B}}
 \int_{\omega'}^{\infty}
 \left\| \Delta H ({\rm i}\omega) \right\|_{\rm F}^2
 \; {\rm d}\omega = 0 . $$
Hence we obtain a frequency $\omega'(\varepsilon)>0$ such that
$$ \max_{\beta \in \mathcal{B}}
 \int_{\omega'(\varepsilon)}^{\infty}
 \left\| \Delta H ({\rm i}\omega) \right\|_{\rm F}^2
 \; {\rm d}\omega < \frac{2\pi}{3} \varepsilon^2 . $$
 The integration domain $(-\infty,-\omega'(\varepsilon))$ exhibits
 the same bound due to a symmetry.

Now we apply Lemma~\ref{lemma:error} to the interval
$[-\omega'(\varepsilon),\omega'(\varepsilon)]$ and obtain
a bound~(\ref{error-compact}) with constants
$K_{\omega'(\varepsilon)},L_{\omega'(\varepsilon)} > 0$.
Let $m = \min \{ n_{\rm in} , n_{\rm out} \}$.
The matrix norms exhibit the general bound
$\| \Delta H \|_F \le \sqrt{m} \| \Delta H \|_2$.
We obtain
$$ \int_{-\omega'(\varepsilon)}^{\omega'(\varepsilon)}
\left\| \Delta H ({\rm i}\omega) \right\|_{\rm F}^2 \; {\rm d}\omega <
2 m \omega'(\varepsilon)
K_{\omega'(\varepsilon)}^2 \beta^2 < \frac{2\pi}{3} \varepsilon^2 $$
provided that $\beta < \Xi(\varepsilon)$ with the constant
$$ \Xi(\varepsilon) = \min \left\{
\frac{\sqrt{\pi} \, \varepsilon}{\sqrt{3 m \omega'(\varepsilon)}
  K_{\omega'(\varepsilon)}} ,
L_{\omega'(\varepsilon)} \right\} . $$
Thus the estimate~(\ref{error-global}) is satisfied for all
$0 < \beta < L$ with the single constant
$L = \min\{ \Xi(\varepsilon), \beta' \}$
and $\beta' < \min \{ \rho(A')^{-2} , 1 \}$.
\hfill $\Box$
%%%%%%%%%%%%%%%%%%%%%%%%%%%%%%%%%%%%%%%%%%%%%%%%%%%%%%%%%%%%%%%%%%%%%%%%%%%%%%%

\medskip

Theorem~\ref{thm:error} illustrates that a regularisation with a 
low error in the $\htwo$-norm can be achieved provided that the
parameters~$\alpha,\beta$ are chosen sufficiently small.
The derived constants are pessimistic, because rough estimates
appear in the proof.
The magnitude of appropriate regularisation parameters
depends on the system of DAEs~(\ref{linear-system}).
However, tiny parameters cause problems in numerical computations
due to ill-conditioned matrices, for example.

%For the lower frequencies, the bound~(\ref{error-compact})
%of Lemma~\ref{lemma:error} yields
%$$  \int_{-\omega'}^{\omega'}
%\left\| H_{\rm DAE}({\rm i}\omega) - H_{\rm ODE}({\rm i}\omega)
%\right\|_{\rm F}^2
%\; {\rm d}\omega < 2 \omega' K_{\omega'}^2 \max\{ \alpha,\beta \}^2 . $$
If the system of DAEs is just proper, then the polynomial part
in~(\ref{transfer-dae}) represents a (non-zero) constant.
Consequently, both the $\htwo$-error and the $\hinf$-error of the
regularisation do not become arbitrarily small.
There is still some potential for using this regularisation technique.
The input-output relation of the linear dynamical system~(\ref{linear-system})
is given by $Y(s) = H_{\rm DAE}(s) U(s)$ with the Laplace transforms
$U,Y$ of input and output, respectively.
If a particular input induces a Laplace transform with a sufficiently fast
decay for high frequencies, then the same effect as in a strictly proper
system emerges.

%%%%%%%%%%%%%%%%%%%%%%%%%%%%%%%%%%%%%%%%%%%%%%%%%%%%%%%%%%%%%%%%%%%%%%%%%%%%%
%%%                      Numerical Experiments                            %%%
%%%%%%%%%%%%%%%%%%%%%%%%%%%%%%%%%%%%%%%%%%%%%%%%%%%%%%%%%%%%%%%%%%%%%%%%%%%%%

\section{Numerical experiments}
\label{sec:examples}
We apply the stability-preserving technique from
Section~\ref{sec:preservation} to two high-dimen\-sional
examples now.
All numerical computations were executed by the software package
MATLAB~\cite{matlab2018}.

%%%%%%%%%%%%%%%%%%%%%%%%%%%%%%%%%%%%%%%%%%%%%%%%%%%%%%%%%%%%%%%%%%%%%%%%%%%%%
\subsection{Microthruster benchmark}
\label{sec:microthruster}
\begin{figure}
  \begin{center}
  \includegraphics[width=6.5cm]{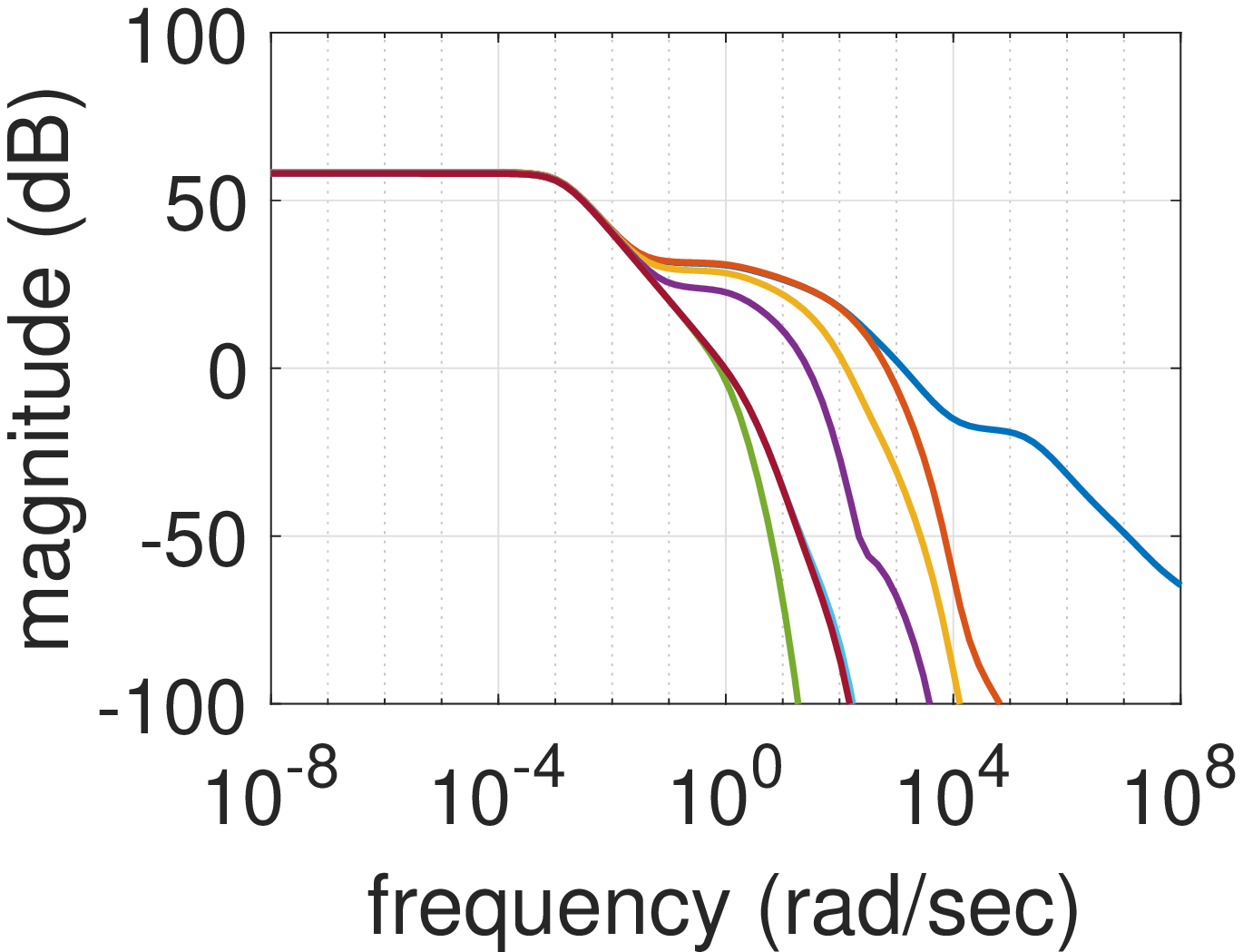}
  \hspace{5mm}
  \includegraphics[width=6.5cm]{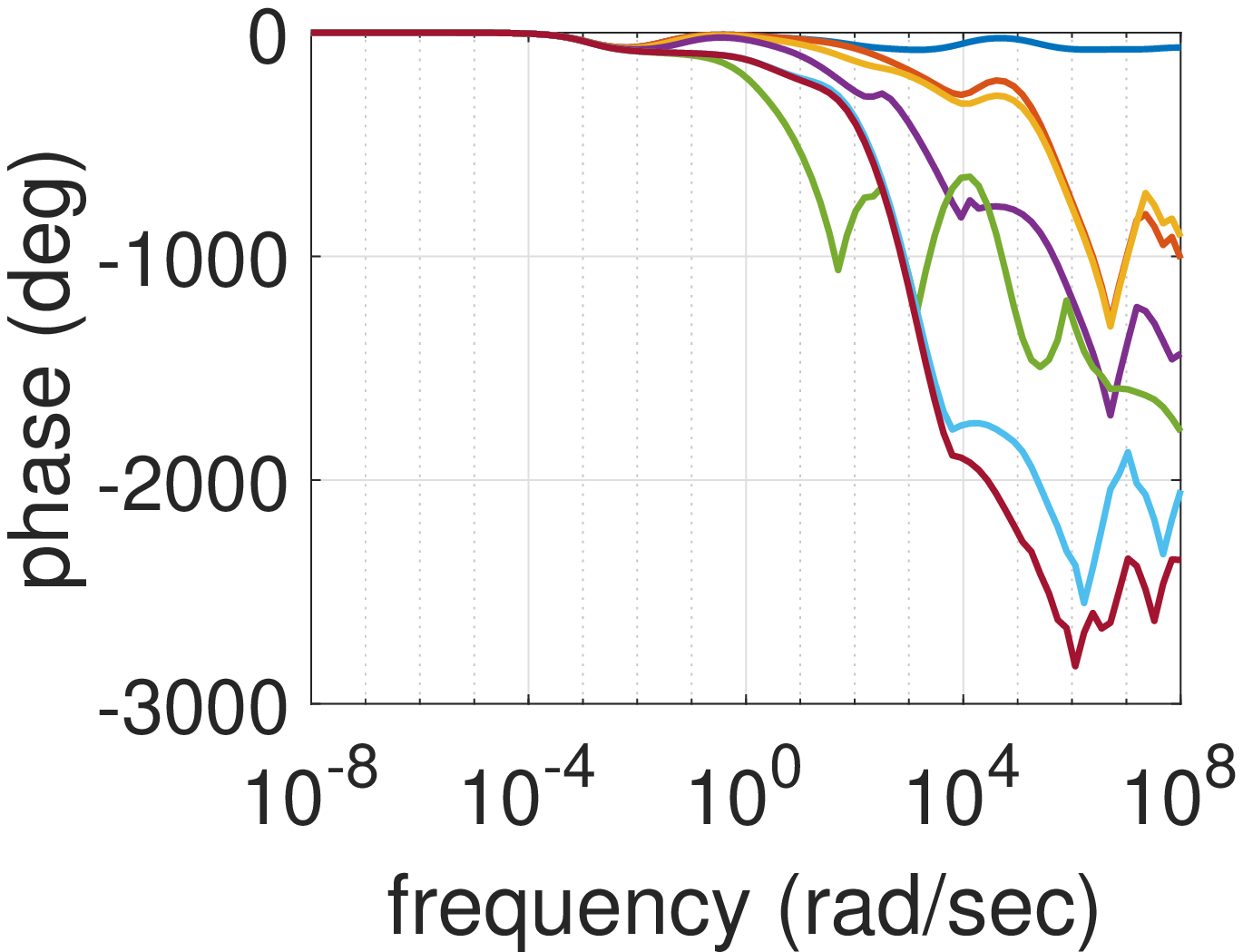}
  \end{center}
  \caption{Bode plot of microthruster example with seven outputs.}
\label{fig:microthruster-bode}
\end{figure}

In~\cite{morwiki}, a microthruster unit represents a test example
called boundary condition independent thermal model.
A spatial discretisation of the two-dimensional
heat transfer partial differential equation
yields a system of ODEs.
More details can be found in~\cite{rudnyi}.
Three parameters appear in this model, which we choose all equal to one.
The linear dynamical system~(\ref{linear-system})
is single-input-multiple-output (SIMO) with $n_{\rm out} = 7$.
Its Bode plot is depicted in Figure~\ref{fig:microthruster-bode}. 
Table~\ref{tab:micro} illustrates some properties of the system.
In particular, the system is asymptotically stable.
The matrix~$E$ is diagonal with positive elements.
Thus we simply scale this system into explicit ODEs ($E=I_n$),
which are used in the following.

%%% Table: Microthruster System %%%%%%%%%%%%%%%%%%%%%%%%%%%%%%%%%%%%%%%%%%%%%
\begin{table}
  \caption{Properties of the microthruster benchmark system.\label{tab:micro}}
\begin{center}
  \begin{tabular}{cc} \hline
    dimension~$n$ & 4257 \\
    \# outputs & 7 \\
    \# non-zero entries in~$A$ & 37465 \\
    \# non-zero entries in~$E$ & 4257 \\
    spectral abscissa $\alpha(E^{-1}A)$ & $-0.0013$ \\ \hline
  \end{tabular}
\end{center}
\end{table}
%%%%%%%%%%%%%%%%%%%%%%%%%%%%%%%%%%%%%%%%%%%%%%%%%%%%%%%%%%%%%%%%%%%%%%%%%%%%%

%%% Figure: Bode Plot %%%%%%%%%%%%%%%%%%%%%%%%%%%%%%%%%%%%%%%%%%%%%%%%%%%%%%%
%%%%%%%%%%%%%%%%%%%%%%%%%%%%%%%%%%%%%%%%%%%%%%%%%%%%%%%%%%%%%%%%%%%%%%%%%%%%%

\begin{figure}
  \begin{center}
  \includegraphics[width=10cm]{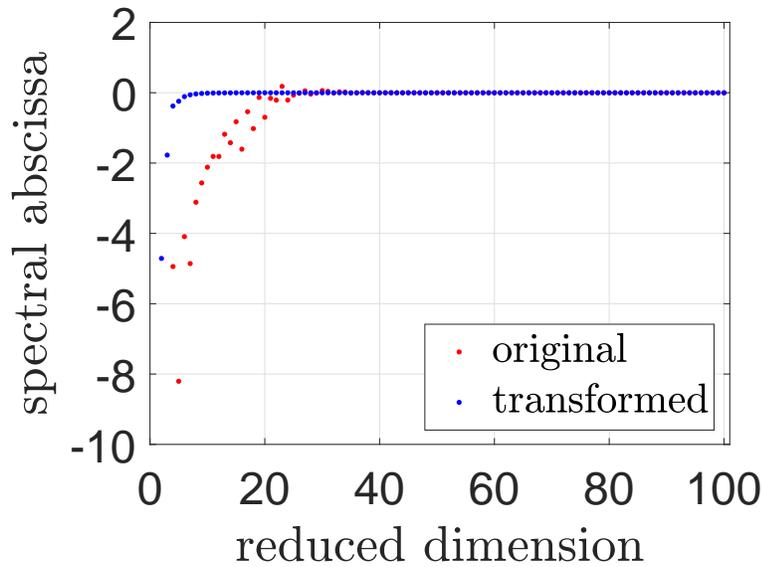}
  \end{center}
  \caption{Spectral abscissa of the matrices in the ROMs
    from original system and transformed system in microthruster example.}
\label{fig:micro-spectrum}
\end{figure}

We use the one-sided Arnoldi method with the single real expansion
point $s_0 = 100$.
The reduced systems are arranged for $r=1,2,\ldots,100$.
The spectral abscissas of the reduced systems are depicted
in Figure~\ref{fig:micro-spectrum}.
It follows that 11 ROMs become unstable, which are all
in the range $20 < r < 50$.
Thus we apply the stability preserving approach from
Section~\ref{sec:preservation}.

%%% Figure: Spectral Abscissa %%%%%%%%%%%%%%%%%%%%%%%%%%%%%%%%%%%%%%%%%%%%%%%
%%%%%%%%%%%%%%%%%%%%%%%%%%%%%%%%%%%%%%%%%%%%%%%%%%%%%%%%%%%%%%%%%%%%%%%%%%%%%

\begin{figure}
  \begin{center}
  \includegraphics[width=6cm]{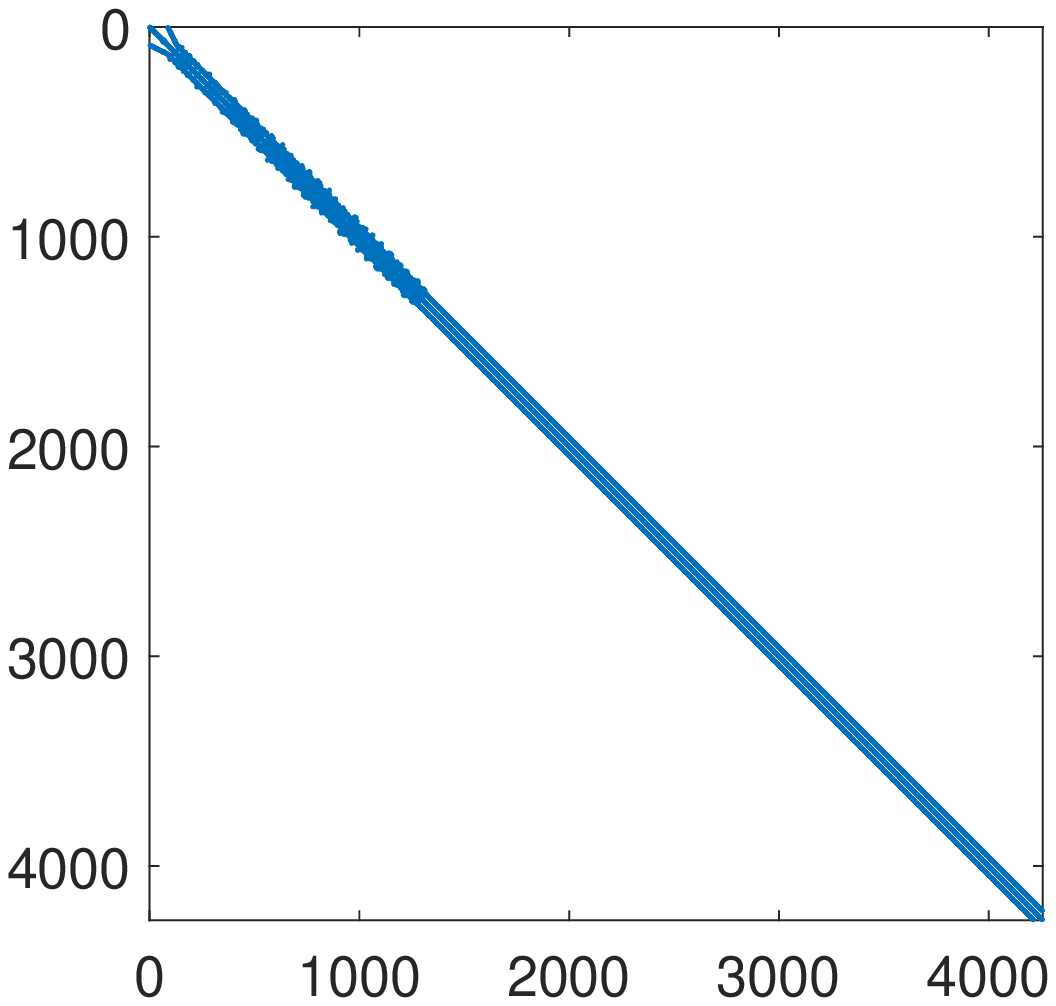}
  \hspace{10mm}
  \includegraphics[width=6cm]{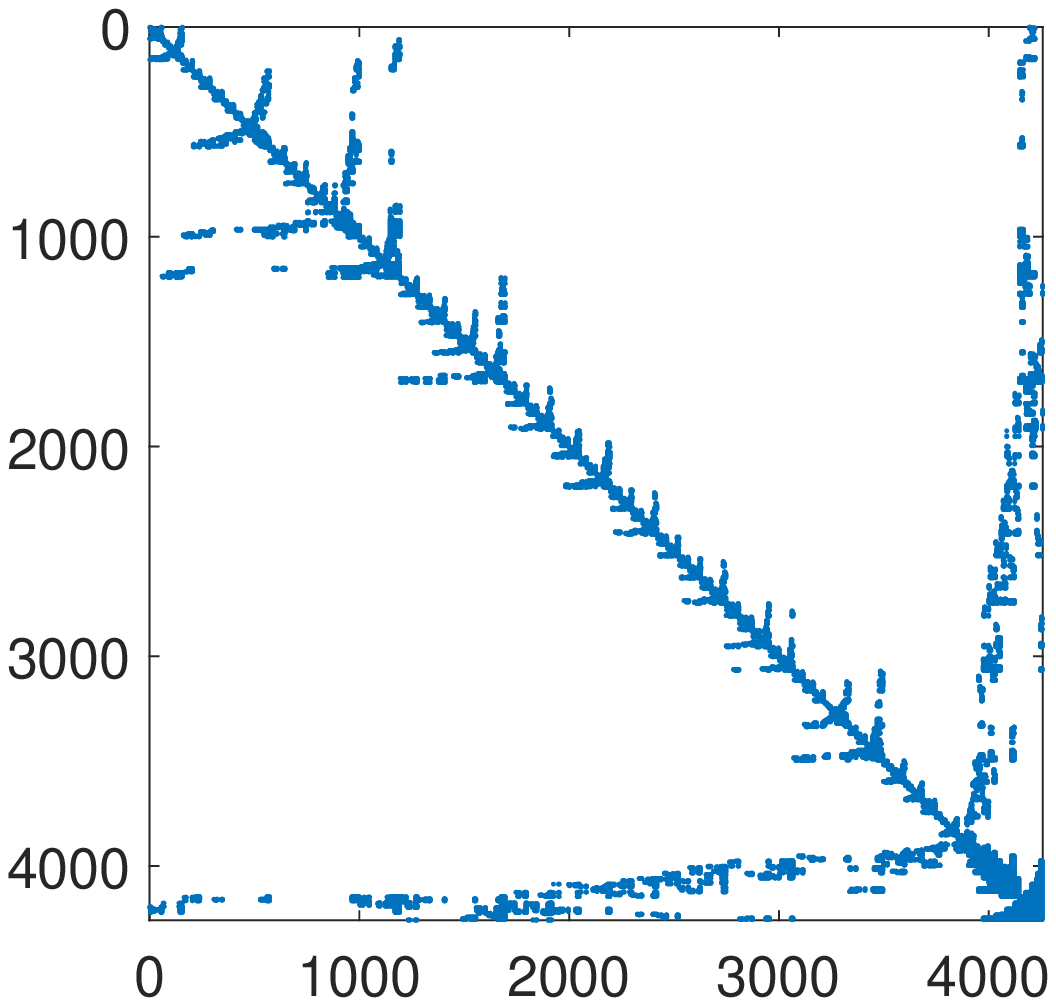}
  \end{center}
  \caption{Sparsity patterns: system matrix~$A$ (left) and 
    $LU$-decomposition of $-{\rm i}I_n-A$ (right).}
\label{fig:micro-sparse}
\end{figure}

A quadrature method requires the solution of
linear systems~(\ref{linear-system-alg}).
We investigate the sparsity in the $LU$-decomposition~(\ref{lu-decomp})
of the matrix
$-{\rm i}I_n-A$ ($\omega=1$).
The number of non-zeros in the $LU$-decomposition with partial pivoting
($Q=I_n$) is $395 \, 463$.
Alternatively, UMFPACK achieves a factorisation with
$226 \, 395$ non-zero entries, whose sparsity pattern
is shown in Figure~\ref{fig:micro-sparse} (right).

%%% Figure: Sparsity Pattern %%%%%%%%%%%%%%%%%%%%%%%%%%%%%%%%%%%%%%%%%%%%%%%%
%%%%%%%%%%%%%%%%%%%%%%%%%%%%%%%%%%%%%%%%%%%%%%%%%%%%%%%%%%%%%%%%%%%%%%%%%%%%%

We investigate three quadrature methods for the computation of
the projection matrix~(\ref{approx-W}):
\begin{itemize}
\item[i)]
  adaptive Gauss-7-Kronrod-15 quadrature using the built-in
  MATLAB function {\tt integral}, see~\cite{shampine},
\item[ii)]
  Gauss-Legendre rule, see~\cite[p.~171]{stoerbulirsch},
  with fixed numbers of nodes,
\item[iii)]
  nested midpoint rules as described in Section~\ref{sec:finite-interval}.
\end{itemize}
In the adaptive quadrature (i), we choose the absolute and the relative
error tolerances as $\varepsilon_{\rm abs} = \varepsilon_{\rm rel} = 0.1$.
The algorithm performs $K=150$ evaluations of the integrand.
The projection matrices~(\ref{matrix-tilde-W}) are used to attain
biorthogonality.
All reduced systems become stable now.
Figure~\ref{fig:micro-spectrum} depicts the spectral abscissas of
the ROMs.
In the Gauss-Legendre rule~(ii), we increase the number of nodes~$K$ until
all ROMs become stable, see Table~\ref{tab:micro-gauss}.
Just $K=14$ nodes are sufficient to obtain always stable systems.
Table~\ref{tab:micro-midpoint} shows the number of stable ROMs for
the refinement in the midpoint rule~(iii).
Now $127$ nodes are required to achieve the stability preservation
in all reduced systems. 

%%% Table: Quadrature %%%%%%%%%%%%%%%%%%%%%%%%%%%%%%%%%%%%%%%%%%%%%%%%%%%%%%%
\begin{table}[h]
  \caption{Number of stable reduced systems out of 100 for different
    numbers of nodes in Gauss-Legendre quadrature.\label{tab:micro-gauss}}
\begin{center}
  \begin{tabular}{rcccccccccccccc} 
    \# nodes & 1 & 2 & 4 & 5 & 6 & 7 & 9 & 11 & 14 \\ \hline
    \# stable systems & 91 & 92 & 93 & 94 & 96 & 97 & 98 & 99 & 100
%    \# nodes & 1 & 2 & 3 & 4 & 5 & 6 & 7 & 8 & 9 & 10 & 11 & 12 & 13 & 14 \\ \hline
%    \# stable systems & 91 & 92 & 93 & 93 & 94 & 96 & 97 & 97 & 98 & 98 & 99 & 99 & 99 & 100  
  \end{tabular}
\end{center}
\end{table}
%%%%%%%%%%%%%%%%%%%%%%%%%%%%%%%%%%%%%%%%%%%%%%%%%%%%%%%%%%%%%%%%%%%%%%%%%%%%%

%%% Table: Quadrature %%%%%%%%%%%%%%%%%%%%%%%%%%%%%%%%%%%%%%%%%%%%%%%%%%%%%%%
\begin{table}[h]
  \caption{Number of stable reduced systems out of 100 for
    nested midpoint rule.\label{tab:micro-midpoint}}
\begin{center}
  \begin{tabular}{rccccccc} 
    \# nodes & 1 & 3 & 7 & 15 & 31 & 63 & 127 \\ \hline
    \# stable systems & 91 & 93 & 93 & 95 & 98 & 99 & 100 
  \end{tabular}
\end{center}
\end{table}
%%%%%%%%%%%%%%%%%%%%%%%%%%%%%%%%%%%%%%%%%%%%%%%%%%%%%%%%%%%%%%%%%%%%%%%%%%%%%

\begin{figure}
\begin{center}
\includegraphics[width=10cm]{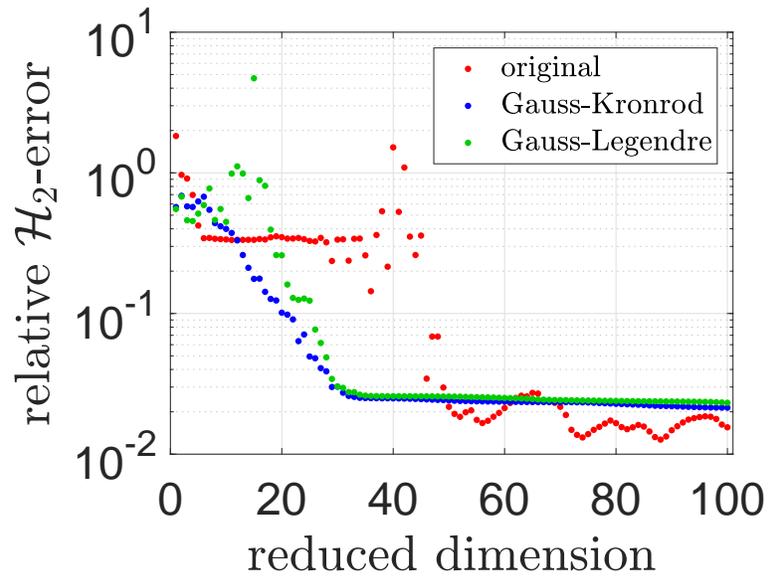} 
\end{center}
\caption{Relative differences in $\htwo$-norm for three MOR approaches
  in microthruster example.}
\label{fig:micro-error}
\end{figure}

We also compare the approximation quality between the ROMs
obtained from the conventional reduction and the stabilisation method.
Figure~\ref{fig:micro-error} illustrates the
relative error in the $\htwo$-norm~(\ref{h2-norm}), i.e.,
\begin{equation} \label{relative-h2-error}
  E_{\rm REL} =
  \frac{\left\| H_{\rm FOM} - H_{\rm ROM} \right\|_{\htwo}}{\left\| H_{\rm FOM} \right\|_{\htwo}}
\end{equation}
including the transfer functions.
We calculate approximations of $\htwo$-norms (\ref{h2-norm}) by
the trapezoidal rule on a logarithmically spaced grid
on the imaginary axis.
The ROMs from the adaptive Gauss-Kronrod rule ($K=150$)
and the Gauss-Legendre rule ($K=14$) are examined.
We observe that the error~(\ref{relative-h2-error}) of the
stabilisation method is often much lower
than the conventional method for dimensions $r<50$,
whereas the errors become close for $r>50$.
Furthermore, the adaptive quadrature is more accurate than the
Gauss-Legendre rule with the low number of nodes for $r<30$.
Our important observation is that the stabilisation approach does
not deteriorate the accuracy of the MOR.

%%% Figure: Error Bound %%%%%%%%%%%%%%%%%%%%%%%%%%%%%%%%%%%%%%%%%%%%%%%%%%%%%
%%%%%%%%%%%%%%%%%%%%%%%%%%%%%%%%%%%%%%%%%%%%%%%%%%%%%%%%%%%%%%%%%%%%%%%%%%%%%

%%%%%%%%%%%%%%%%%%%%%%%%%%%%%%%%%%%%%%%%%%%%%%%%%%%%%%%%%%%%%%%%%%%%%%%%%%%%%
\subsection{Random low-pass filter}
\begin{figure}
  \begin{center}
  \includegraphics[width=14cm]{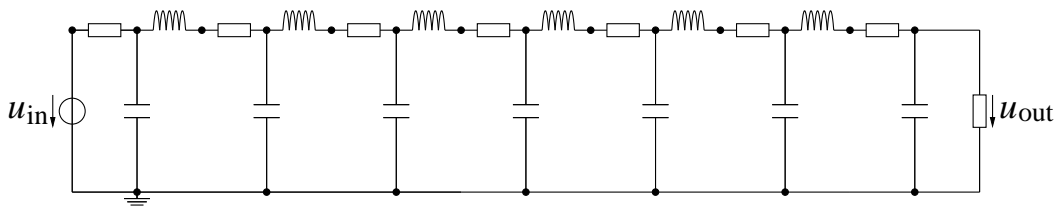}
  \end{center}
  
  \caption{Circuit diagram of low-pass filter.}
  
\label{fig:filter-circuit}
\end{figure}

We investigate the electric circuit of a low-pass filter
in Figure~\ref{fig:filter-circuit}.
The circuit includes 21 physical parameters:
seven capacitances, six inductances and eight
conductances.
A mathematical modelling generates a system of DAEs~(\ref{linear-system})
for 14~node voltages and 6~branch currents ($n'=20$).
The (nilpotency) index of this system is one.
Furthermore, the system is asymptotically stable and strictly proper.
The system is single-input-single-output (SISO),
because a voltage source is supplied and the output is defined
as the voltage at a load conductance.

%%% Figure: Low-Pass-Filter %%%%%%%%%%%%%%%%%%%%%%%%%%%%%%%%%%%%%%%%%%%%%%%%%
%%%%%%%%%%%%%%%%%%%%%%%%%%%%%%%%%%%%%%%%%%%%%%%%%%%%%%%%%%%%%%%%%%%%%%%%%%%%%

\begin{figure}
  \begin{center}
    \includegraphics[width=6.5cm]{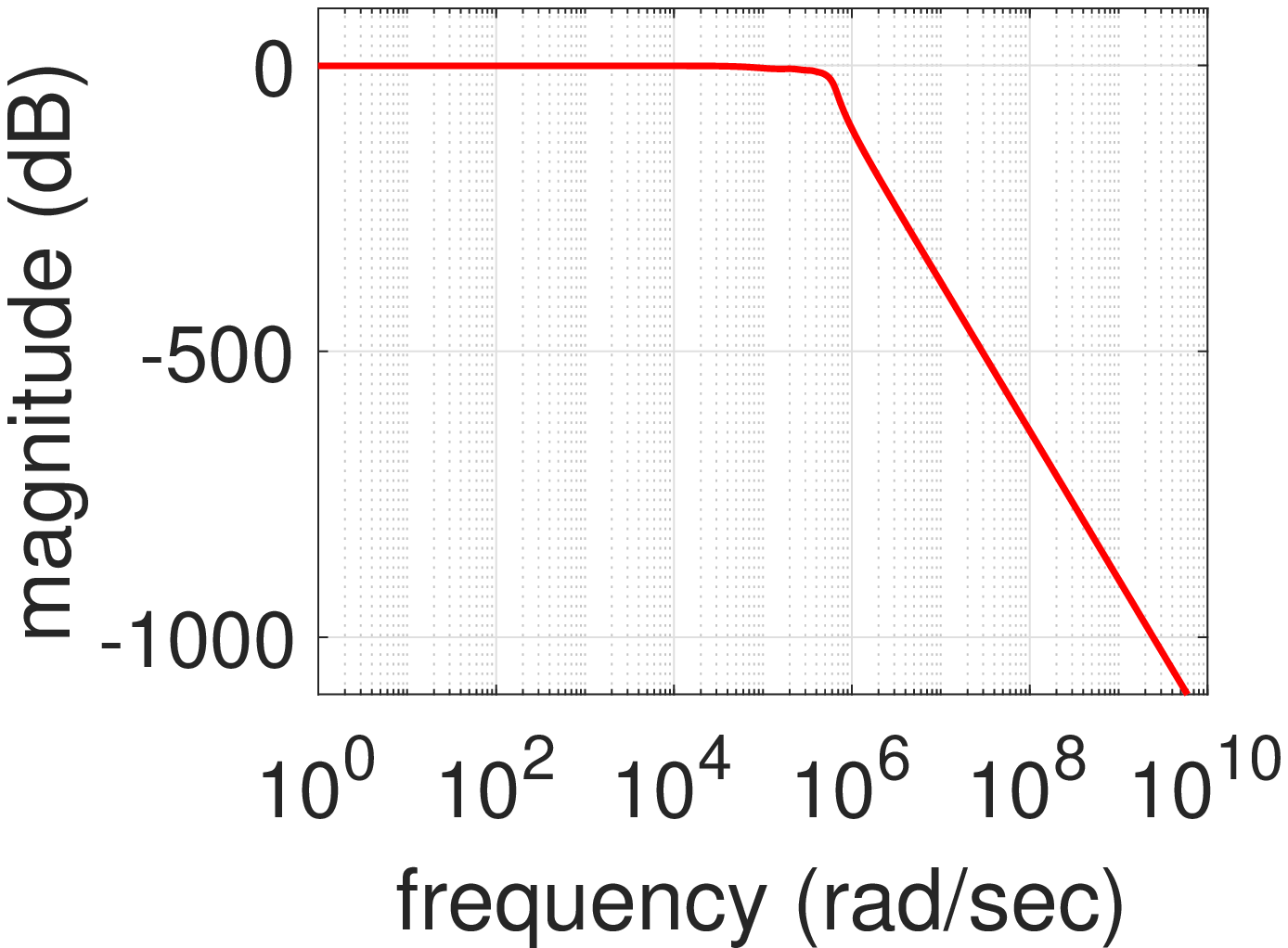}
    \hspace{5mm}
    \includegraphics[width=6.5cm]{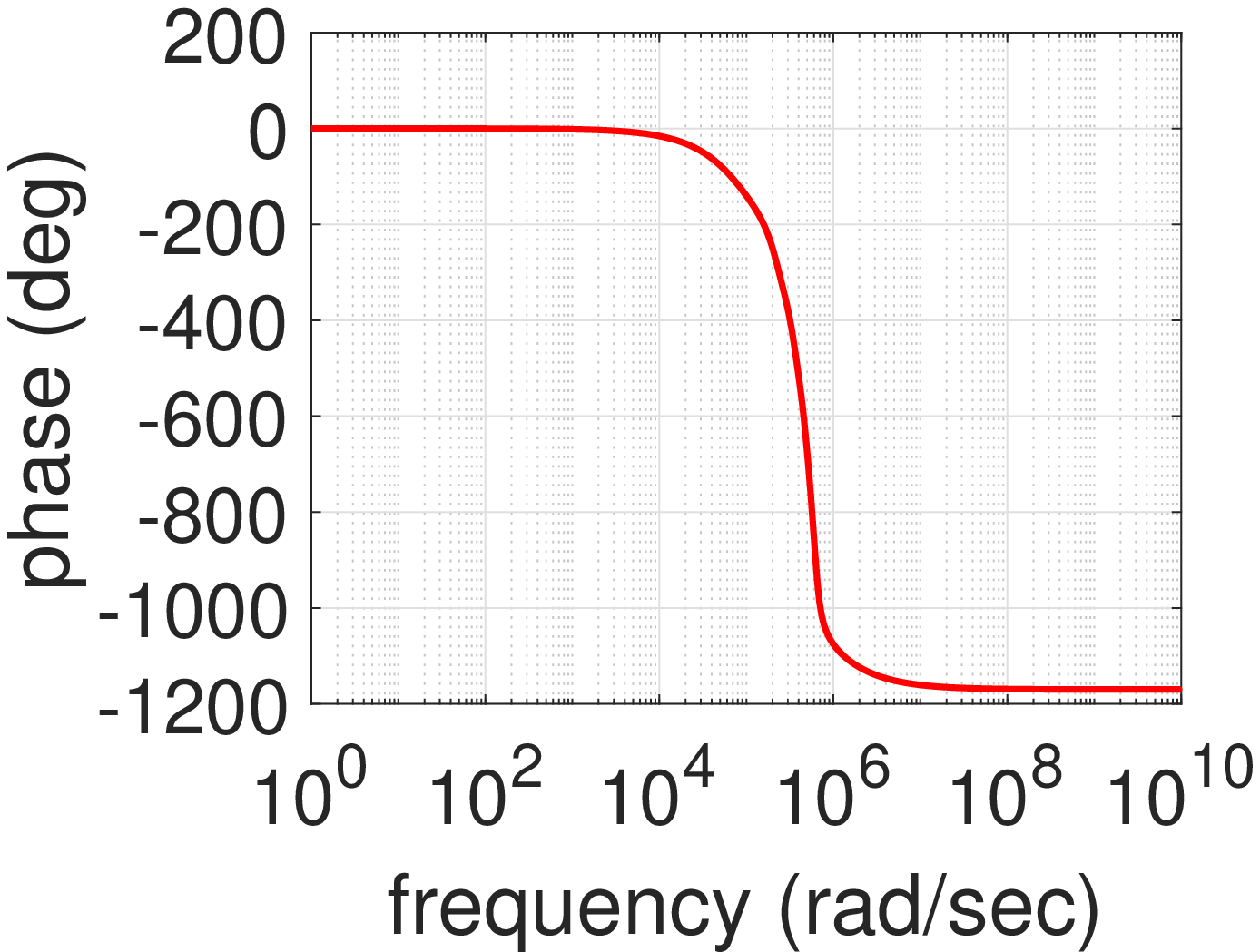}
  \end{center}
  \caption{Bode plot of first output in stochastic Galerkin system
    of random low-pass filter.}
\label{fig:filter-bode}
\end{figure}

In a stochastic modelling, all physical parameters are replaced by
uniformly distributed random variables with a variation of 15\%
around their mean values.
We use a truncated polynomial chaos expansion to approximate the
random processes, see~\cite{xiu-book}.
All basis polynomials are included up to total degree three, i.e.,
$m=2024$ basis functions depending on 21~variables. 
The stochastic Galerkin method yields a larger system of
DAEs~(\ref{linear-system}) with dimension~$n=mn'$, whose
solution approximates the unknown coefficient functions. 
Table~\ref{tab:filter} depicts its characteristic numbers.
The system is SIMO with a large number of outputs. 
This linear dynamical system inherits the properties of the circuit model:
index-one, asymptotically stable and strictly proper.
Figure~\ref{fig:filter-bode} illustrates the Bode plot of
the first output, which represents an approximation for the expected value
of the output voltage.
The magnitude of the transfer function shows that high frequencies
are damped out.
This test example was also used in~\cite{pulch18}.

%%% Table: Low-Pass-Filter %%%%%%%%%%%%%%%%%%%%%%%%%%%%%%%%%%%%%%%%%%%%%%%%%%
\begin{table}
  \caption{Properties of stochastic Galerkin system for
    random low-pass filter. \label{tab:filter}}
\begin{center}
  \begin{tabular}{cc} \hline
    dimension~$n$ & 40480 \\
    \# outputs & 2024 \\
    \# non-zero entries in~$A$ & 116886 \\
    \# non-zero entries in~$E$ & 32890 \\
    rank($E$) & 26312 \\ \hline
  \end{tabular}
\end{center}
\end{table}
%%%%%%%%%%%%%%%%%%%%%%%%%%%%%%%%%%%%%%%%%%%%%%%%%%%%%%%%%%%%%%%%%%%%%%%%%%%%%

%%% Figure: Bode Plot %%%%%%%%%%%%%%%%%%%%%%%%%%%%%%%%%%%%%%%%%%%%%%%%%%%%%%%
%%%%%%%%%%%%%%%%%%%%%%%%%%%%%%%%%%%%%%%%%%%%%%%%%%%%%%%%%%%%%%%%%%%%%%%%%%%%%

\begin{figure}
  \begin{center}
    \includegraphics[width=6cm]{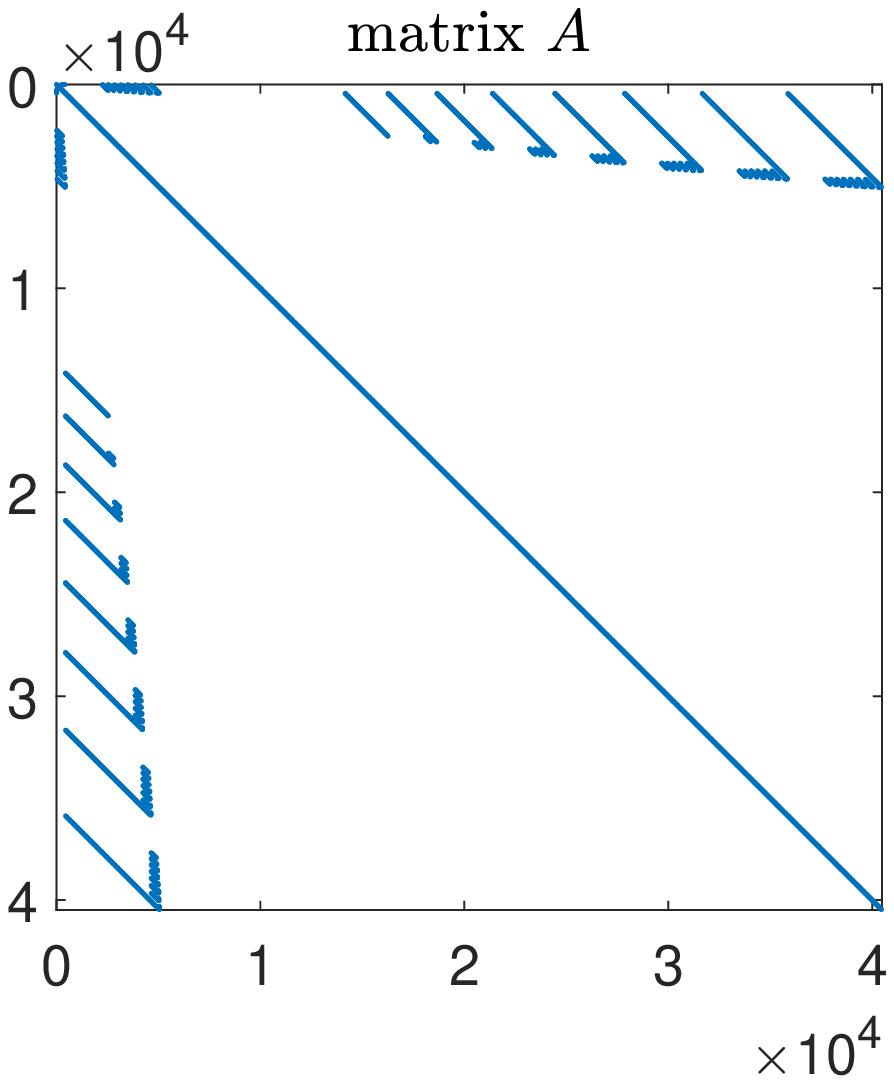}
    \hspace{18mm}
    \includegraphics[width=6cm]{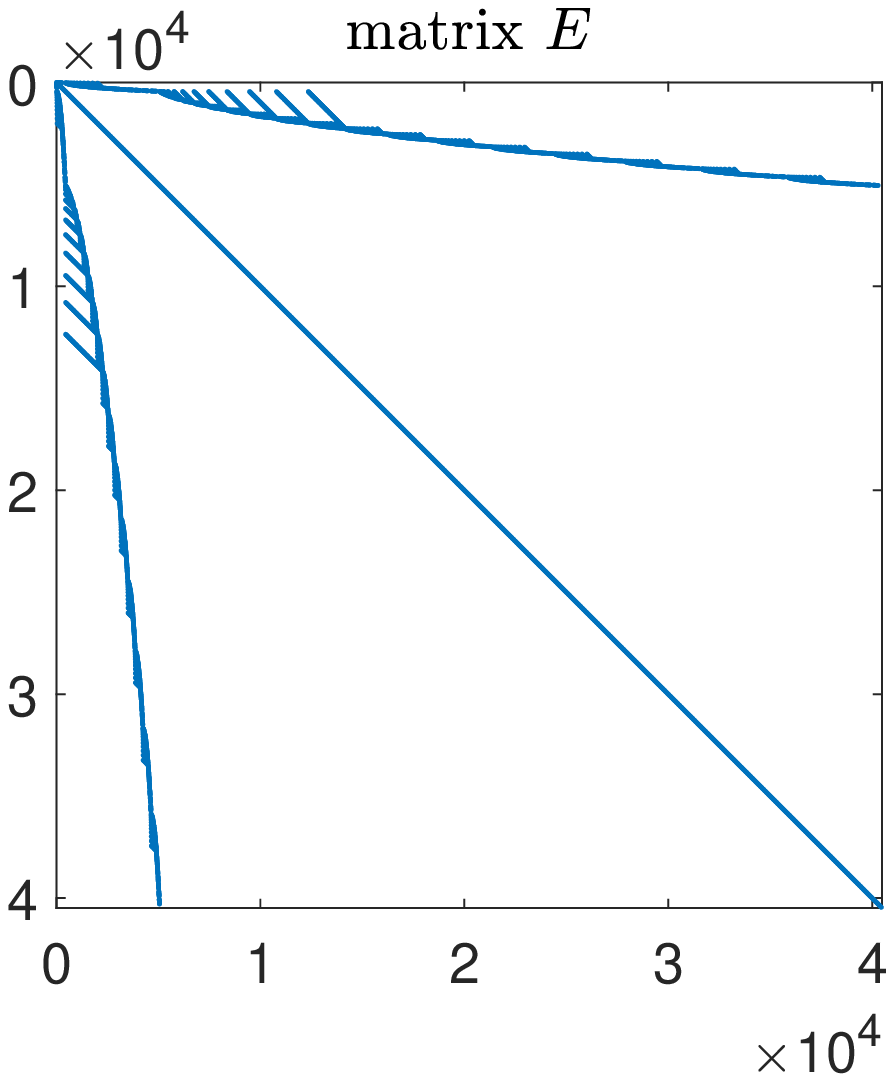}

    \vspace{3mm}

    \includegraphics[width=6cm]{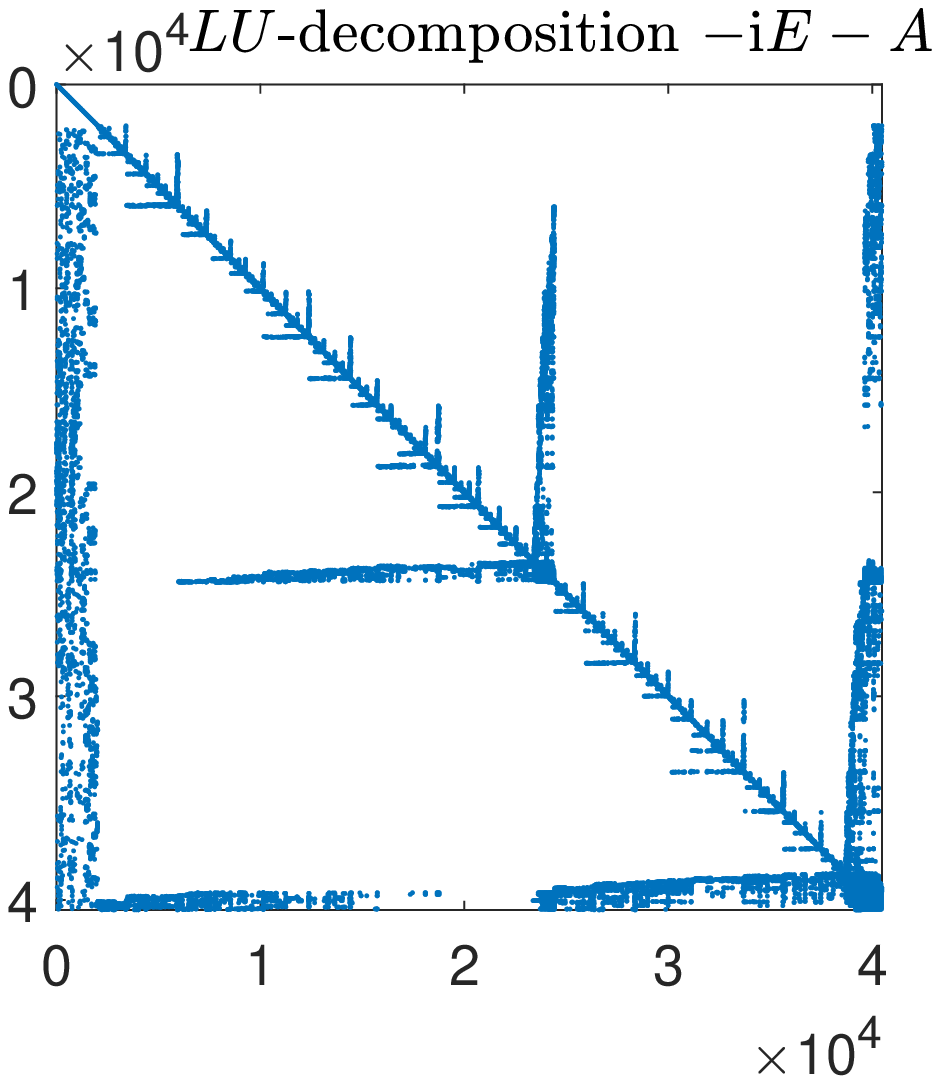}
  \end{center}
  \caption{Sparsity patterns of matrices in stochastic Galerkin system
    of random low-pass filter.}
\label{fig:filter-sparse}
\end{figure}

Furthermore,
Figure~\ref{fig:filter-sparse} depicts the sparsity patterns of the
system matrices and the $LU$-de\-composi\-tion~(\ref{lu-decomp})
for $\omega=1$.
UMFPACK yields an $LU$-factorisation with about $1.8$ million non-zero
entries.
In contrast, the common $LU$-factorisation with partial pivoting
generates about 150 million non-zero elements and thus requires much
more computational work.
Consequently, we employ decompositions from UMFPACK whenever
linear systems of this type appear.

%%% Figure: Sparsity Pattern %%%%%%%%%%%%%%%%%%%%%%%%%%%%%%%%%%%%%%%%%%%%%%%%
%%%%%%%%%%%%%%%%%%%%%%%%%%%%%%%%%%%%%%%%%%%%%%%%%%%%%%%%%%%%%%%%%%%%%%%%%%%%%

\begin{figure}
  \begin{center}
  \includegraphics[width=6.5cm]{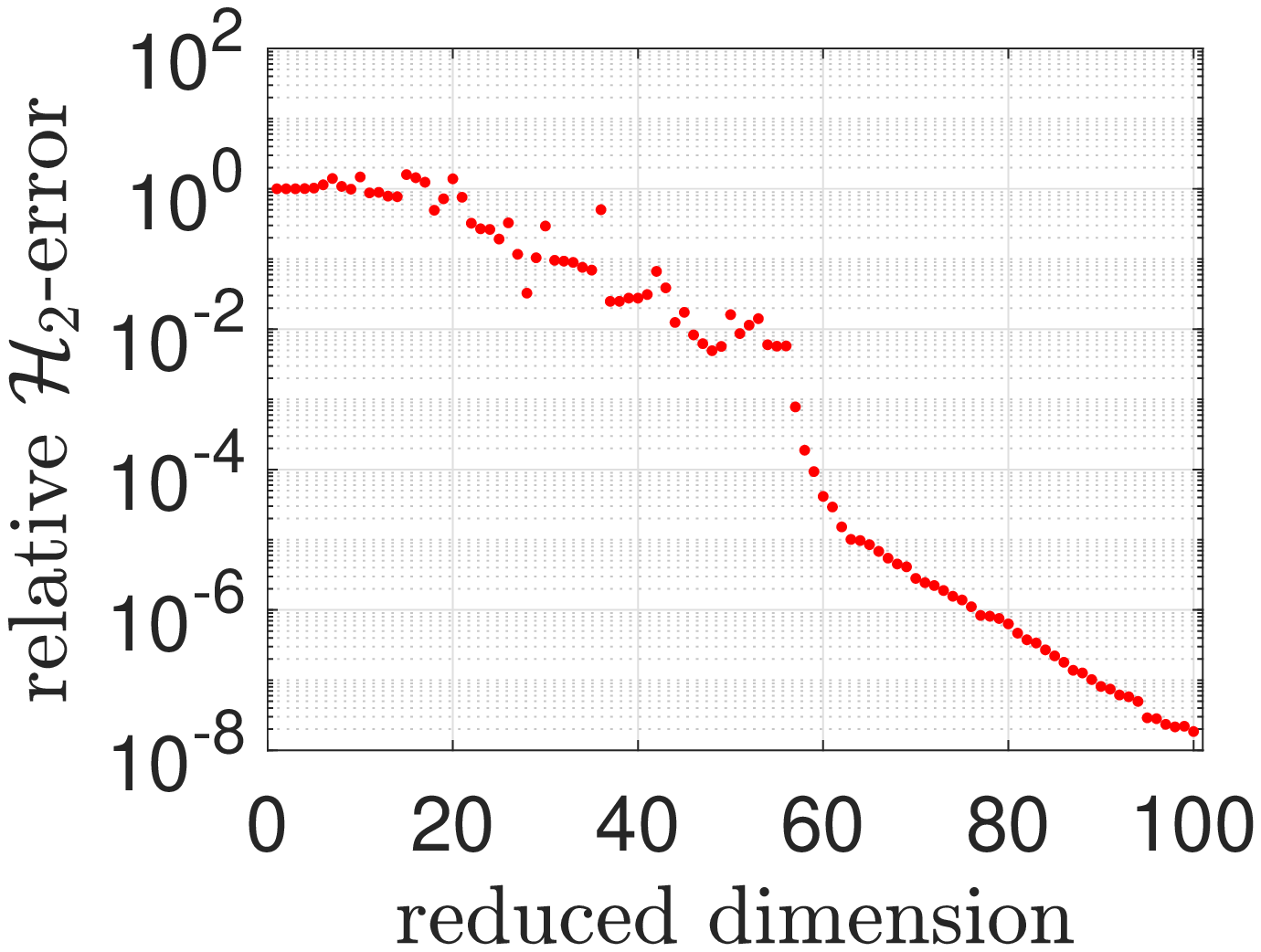}
  \hspace{5mm}
  \includegraphics[width=6.5cm]{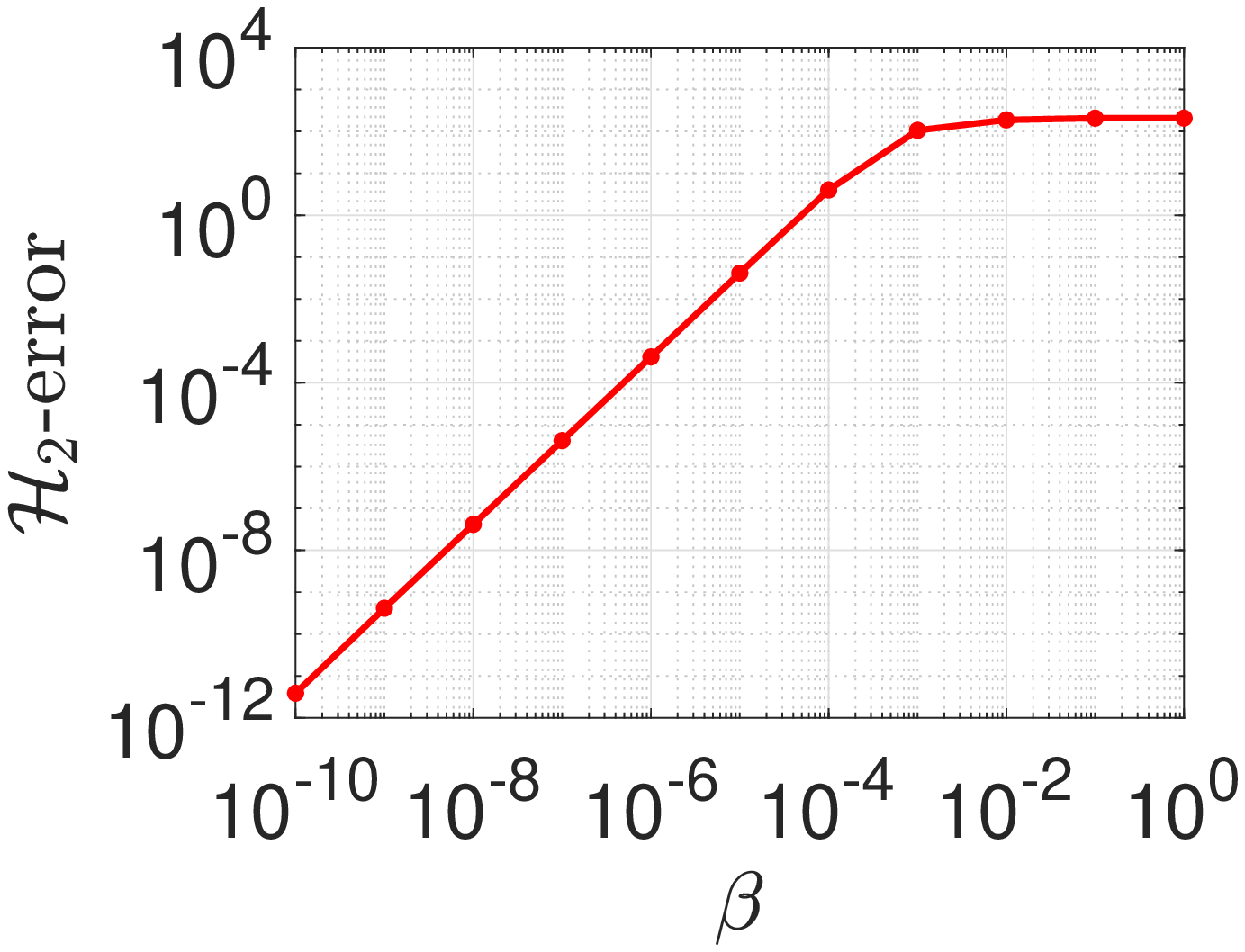}
  \end{center}
  \caption{Error of MOR for DAE system in relative $\htwo$-norm (left)
    and error of the regularisation in $\htwo$-norm (right)
    for stochastic Galerkin system.}
\label{fig:galerkin-errors}
\end{figure}

Now the one-sided Arnoldi method with the single real expansion point
$s_0 = 5 \cdot 10^5$ is used in all cases.
The ROMs are always computed for dimensions $r=1,2,\ldots,100$.
We directly reduce the system of DAEs first.
The relative $\htwo$-errors~(\ref{relative-h2-error}) of the MOR
are depicted in Figure~\ref{fig:galerkin-errors} (left).
This error decays rapidly and becomes tiny.
However, only 6 out of 100 ROMs inherit the asymptotic stability
of the FOM.
Hence a stability-preserving method is essential in this example.

The system of DAEs is regularised by the technique described in
Section~\ref{sec:regularisation}.
We arrange the modified matrices~(\ref{matrices-regularised})
with $\alpha = \beta^2$ for different parameters~$\beta$.
The total error~(\ref{dae-total-error}) is bounded by the sum of 
regularisation error and MOR error.
Figure~\ref{fig:galerkin-errors} (right) shows the (absolute)
$\htwo$-error of the regularisation.
We recognise that this error converges exponentially to zero
for $\beta$ tending to zero.
In addition, numerical computations confirm that the investigated
regularised systems are asymptotically stable.

On the one hand, we apply the conventional Arnoldi method to the systems of
ODEs for several parameters~$\beta$.
On the other hand, we perform the stabilisation technique of
Section~\ref{sec:preservation} for the ODEs
in combination with the Arnoldi algorithm.
The adaptive Gauss-Kronrod quadrature yields the associated
projection matrix as in Section~\ref{sec:microthruster}.
The used error tolerances read as
$\varepsilon_{\rm abs}=\varepsilon_{\rm rel}=0.1$ again.
Table~\ref{tab:galerkin-stab} illustrates the number of stable ROMs
and the number of nodes in the quadrature.
The conventional approach generates more and more unstable systems for
decreasing parameters~$\beta$.
In contrast, the stabilised method always yields at least 95\%
stable reduced systems.
Moreover, the unstable ROMs occur only within dimensions $r < 10$.
The number of nodes, which are selected by the adaptive quadrature,
increases for decaying parameters~$\beta \rightarrow 0$.
This behaviour reflects that the 
integral~(\ref{M-times-V}) does not exist in the limit case~$\beta=0$.
However, the ratios $\beta / K$ of the regularisation parameter and
the number of nodes still converges nearly exponentially to zero in the
observed range.
Thus the rise in~$K$ is much lower than the decay in~$\beta$.
It turns out that alleviated tolerances
$\varepsilon_{\rm abs},\varepsilon_{\rm rel}$
cause more unstable systems.
Furthermore, the Gauss-Legendre rule performs worse in this example.

%%% Table: Stability %%%%%%%%%%%%%%%%%%%%%%%%%%%%%%%%%%%%%%%%%%%%%%%%%%%%%%%%
\begin{table}
  
  \caption{Number of stable ROMs for different regularisation parameters.}
\begin{center}
  \begin{tabular}{rcccccc}
    parameter~$\beta$ & $10^{-2}$ & $10^{-3}$ & $10^{-4}$ & $10^{-5}$ & $10^{-6}$ & $10^{-7}$ \\ \hline
    \# stable, conventional & 100 & 77 & 56 & 57 & 51 & 43 \\
    \# stable, stabilised & 100 & 99 & 95 & 95 & 95 & 95 \\
    \# nodes in quadrature & 180 & 330 & 480 & 600 & 810 & 900 
  \end{tabular}
\end{center}

\label{tab:galerkin-stab}
\end{table}
%%%%%%%%%%%%%%%%%%%%%%%%%%%%%%%%%%%%%%%%%%%%%%%%%%%%%%%%%%%%%%%%%%%%%%%%%%%%%

%%% Figure: Errors in stoch. Galerkin method  %%%%%%%%%%%%%%%%%%%%%%%%%%%%%%%
%%%%%%%%%%%%%%%%%%%%%%%%%%%%%%%%%%%%%%%%%%%%%%%%%%%%%%%%%%%%%%%%%%%%%%%%%%%%%

\begin{figure}
  \begin{center}
  \includegraphics[width=6.5cm]{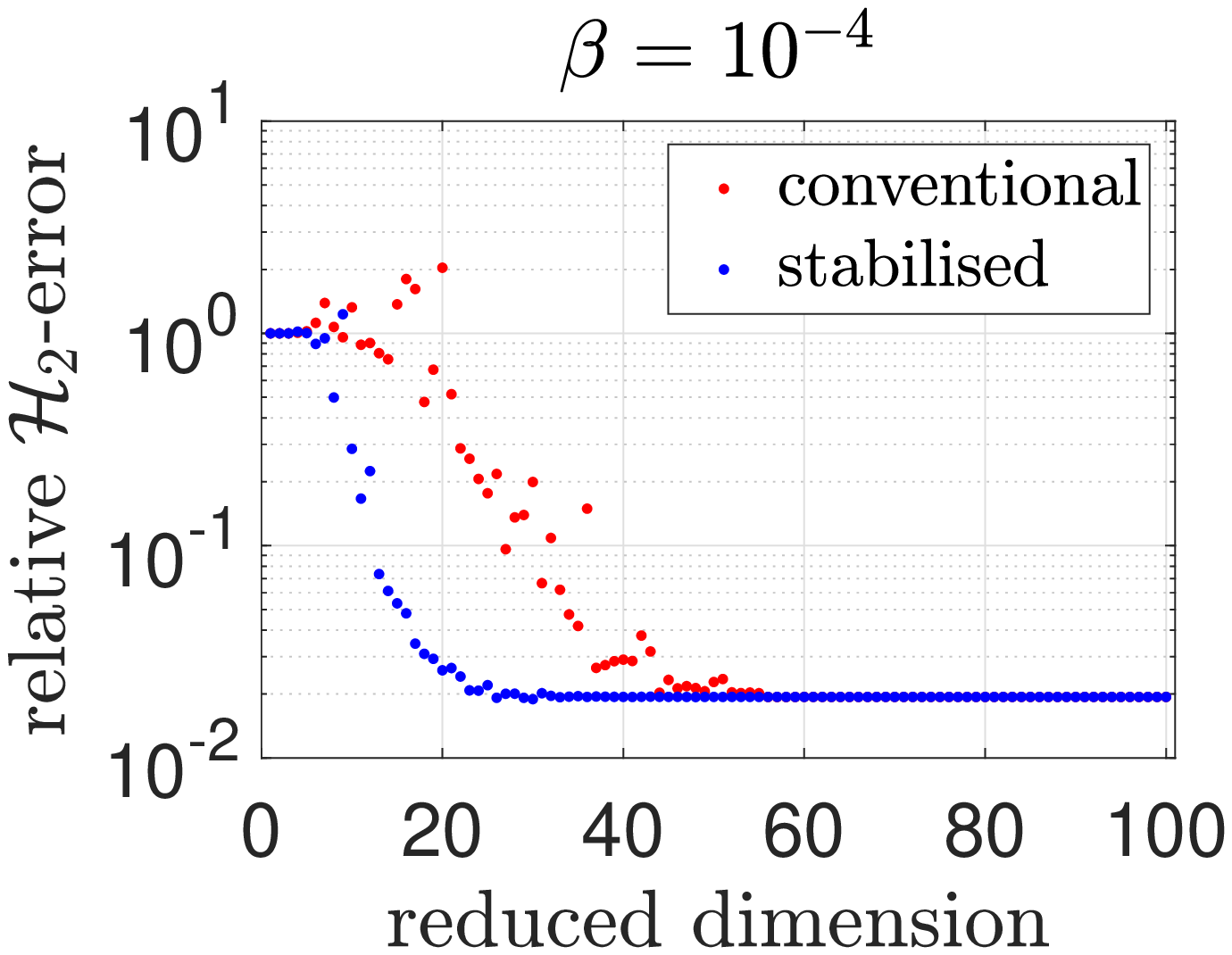}
  \hspace{5mm}
  \includegraphics[width=6.5cm]{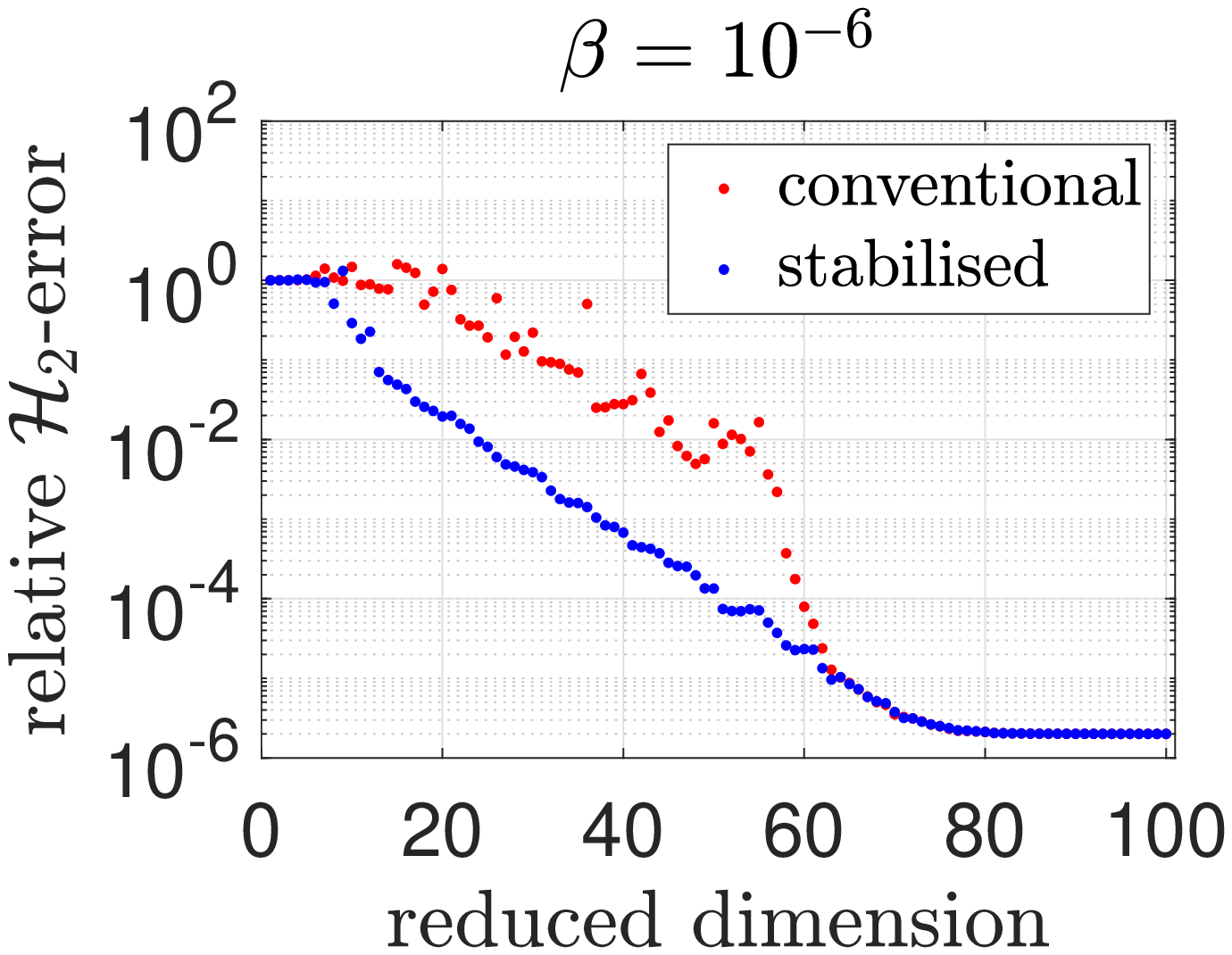}
  \end{center}
  \caption{Error in MOR of regularised systems from stochastic Galerkin
    method for parameters
    $\beta = 10^{-4}$ (left) and $\beta = 10^{-6}$ (right).}
\label{fig:galerkin-mor-ode-error}
\end{figure}

Finally, we examine the total error~(\ref{dae-total-error}) of the MOR,
where the system of DAEs represents the FOM.
Figure~\ref{fig:galerkin-mor-ode-error} shows the relative errors
with respect to the $\htwo$-norm in the two cases
$\beta = 10^{-4},10^{-6}$.
The error decreases fastly for low reduced dimensions.
Thereafter the total error stagnates,
because the error of the regularisation dominates.
We observe that the total error of the stability-preserving approach
is always smaller or equal in comparison to the conventional technique.
Moreover, the stabilised MOR method yields much smaller errors
in the case of low dimensions.

%%% Figure: Error Bound %%%%%%%%%%%%%%%%%%%%%%%%%%%%%%%%%%%%%%%%%%%%%%%%%%%%%
%%%%%%%%%%%%%%%%%%%%%%%%%%%%%%%%%%%%%%%%%%%%%%%%%%%%%%%%%%%%%%%%%%%%%%%%%%%%%

\clearpage

%%%%%%%%%%%%%%%%%%%%%%%%%%%%%%%%%%%%%%%%%%%%%%%%%%%%%%%%%%%%%%%%%%%%%%%%%%%%%
%%%                          Conclusions                                  %%%
%%%%%%%%%%%%%%%%%%%%%%%%%%%%%%%%%%%%%%%%%%%%%%%%%%%%%%%%%%%%%%%%%%%%%%%%%%%%%

\section{Conclusions}
In Galerkin-type projection-based MOR, 
stability preservation can be achieved by a transformation of
a projection matrix.
The transformation is associated with a high-dimensional Lyapunov inequality,
which is satisfied by solving a specific Lyapunov equation.
We designed a numerical method to compute the alternative projection matrix,
where quadrature methods determine approximations of integrals
in the frequency domain.
In contrast to other numerical solvers of high-dimensional
Lyapunov equations,
our frequency domain integral ensures that the inherent approximate solution
is a non-singular matrix.
Results of numerical computations demonstrate that our approach is
efficient in the case of ODEs.

We also generalised the numerical technique to DAEs by a regularisation.
Again quadrature rules, which are applied to frequency domain integrals, 
yield the projection matrices for the regularised ODEs.
However, the quadrature methods require more and more nodes
for decreasing regularisation parameters,
This property restricts the efficiency of our approach to some extend
in the case of DAEs.

%%%%%%%%%%%%%%%%%%%%%%%%%%%%%%%%%%%%%%%%%%%%%%%%%%%%%%%%%%%%%%%%%%%%%%%%%%%%%
%%%                           References                                  %%%
%%%%%%%%%%%%%%%%%%%%%%%%%%%%%%%%%%%%%%%%%%%%%%%%%%%%%%%%%%%%%%%%%%%%%%%%%%%%%

\end{document}